\documentclass[11pt]{article}
\usepackage{amssymb}
\font\tenBbb=msbm10
\def\Z{\hbox{\tenBbb Z}}
\def\T{{\hbox{\tenBbb T}}}
\font\tenBbb=msbm10

\def\sgn{\hbox{sgn}}
\font\tenBbb=msbm10
\def\N{\hbox{\tenBbb N}}
\font\tenBbb=msbm10
\newtheorem{theorem}{Theorem}[section]
\newtheorem{pro}{Proposition}[section]
\newtheorem{lem}{Lemma}[section]

\newtheorem{rem}{Remark}[section]
\newcommand{\R}{{ I\!\!R}}

\def\supp{\mathop{\rm supp}\nolimits}

\newcommand{\re}[1]{(\ref{#1})}

\def\into{\int \hspace*{-4mm} - \,}
\begin{document}
\begin{center}
\noindent {\Large \bf{Global well-posedness in the Energy space
 for  the  Benjamin-Ono equation on the circle}}
\end{center}
\vskip0.2cm
\begin{center}
\noindent
{\bf Luc Molinet}\\
{\small L.A.G.A., Institut Galil\'ee, Universit\'e Paris-Nord,\\
93430 Villetaneuse, France.} \vskip0.3cm
E-mail : molinet@math.univ-paris13.fr.
\end{center}
\vskip0.5cm \noindent {\bf Abstract.} {\small We prove that the Benjamin-Ono equation is locally well-posed in $ H^{1/2}(\T) $. This leads to a global well-posedeness result in $ H^{1/2}(\T) $ thanks to the energy conservation.} \vspace{2mm}\\

\section{Introduction, main results and notations}
\subsection{Introduction}
This paper is devoted to the study of the Cauchy problem
for the
 Benjamin-Ono equation on the circle
$$
 \left\{\begin{array}{llll}
\partial_tu +{\cal{H}}\partial^2_x u - u \partial_x u=0 \,  \, ,\; (t,x)\in \R\times \T \;,\\
u(0,x)=u_0(x) \;,
\end{array}\right. \leqno{\mbox{(BO)}}
$$
where  $ \T=\R/2\pi \Z $, $ u $ is real-valued  and  ${\cal{H}}$ is the Hilbert transform defined for
  $ 2\pi  $-periodic functions with mean value zero by
$$
\widehat{{\cal{H}}(f)}(0)= 0 \quad \mbox{ and }\quad \widehat{{\cal{H}}(f)}(\xi)= -i \,\sgn(\xi) \hat{f}(\xi), \quad
\xi\in \Z^* \quad \quad .
$$
\vskip0.3cm

The Benjamin-Ono equation  arises as a model for long internal
gravity waves in deep stratified fluids, see \cite{B}.  This
equation is formally completely integrable (cf. \cite{AF}, \cite{CW}) and thus possesses
an infinite number of conservation laws. These conservation laws
permit to control  the $ H^{n/2} $-norms, $ n\in\N $, and thus to derive
 global well-posedness results in Sobolev spaces. The Cauchy
problem on the real line has been extensively studied these last
years (cf. \cite{Saut}, \cite{ABFS},  \cite{Io}, \cite{P},
\cite{MST}, \cite{KT1}, \cite{KK}). Recently, T. Tao \cite{T} has
pushed the well-posedness theory to $ H^1(\R) $ by using an
appropriate gauge transform. In the periodic setting, the local
well-posedness of (BO) is known in $ H^s(\T) $ for $s>3/2 $ (cf.
\cite{ABFS}, \cite{Io}), by standard compactness methods which do
not take advantage of the dispersive effects of the equation.
Thanks to the conservation laws  mentioned above and an
interpolation  argument, this leads to global well-posedness in $
H^{s}(\T) $ for $ s>3/2$ (cf. \cite{ABFS}). Very recently, F.
Ribaud and the author \cite{MR3} have improved the global
well-posedness result to $ H^1(\T) $ by using the gauge transform
introduced by T. Tao \cite{T} combining with Strichartz estimates
derived in \cite{B} for the Schr\"odinger group on the
one-dimensional torus.\\
The aim of this paper is to improve the local and global well-posedness
to $ H^{1/2}(\T) $ which is the energy space for (BO).  Recall that
the momentum and the energy for the Benjamin-Ono equation are
given by
\begin{equation}
M(u)=\int_{\T}  u^2 \quad \mbox{ and } \quad
E(u)=\frac{1}{2} \int_{\T} |D_x^{1/2} u |^2 +\frac{1}{6}
 \int_{\T}  u^3 \quad . \label{energy}
\end{equation}
Let us  underline that in order to prove qualitative properties as
stability of travelling waves, a well-posedness result  in the
energy space is often very useful. Our strategy is to combine the gauge
transform of T. Tao with estimates in Bourgain spaces.
\subsection{Notations}
For $x,y\in \R$,  $x\sim y$ means that there exists $C_1$,
$C_2>0$ such that $C_1 |x| \leq |y| \leq C_2|x|$ and $x\lesssim y$
means that there exists $C_2>0$ such that $ |x| \leq C_2|y|$. For
a Banach space $ X $, we denote by $ \| \cdot \|_X $ the norm in
$ X $.\\
We will use the same notations as in \cite{CKSTT1} and \cite{CKSTT2} to deal with Fourier transform of space periodic functions with a large period $\lambda $.
 $ (d\xi)_\lambda $ will be the renormalized counting measure on $ \lambda^{-1} \Z $ :
$$
\int a(\xi)\, d\xi = \frac{1}{\lambda} \sum_{\xi\in
\lambda^{-1}\Z } a(\xi) \quad .
$$
As written in \cite{CKSTT2}, $ (d\xi)_\lambda $ is the counting measure on the integers when $ \lambda=1 $ and converges weakly
 to the Lebesgue measure when $ \lambda\to \infty $. In all the text, all the Lebesgue norms in $ \xi $ will be with respect to the
 measure $ (d\xi)_\lambda $. For  a $ (2\pi \lambda) $-periodic function $ \varphi$, we define its space Fourier transform on
  $ \lambda^{-1}\Z$ by
$$
\hat{\varphi}(\xi)=\int_{\R/(2\pi \lambda)\Z} e^{-i \xi x} \, f(x)
\, dx , \quad \forall \xi \in \lambda^{-1}\Z \quad .
$$
We denote by $ V(\cdot) $ the free group associated with the linearized Benjamin-Ono equation,
$$
\widehat{V(t)\varphi}(\xi)=e^{-i \xi|\xi|t} \,
\hat{\varphi}(\xi) , \quad \xi\in \lambda^{-1}\Z \quad .
$$
We define the Sobolev spaces $ H^s_{\lambda} $ for $ (2\pi\lambda)$-periodic functions by
$$
\|\varphi\|_{H^s_\lambda}=\|\langle \xi \rangle^{s} \varphi(\xi)
\|_{L^2_\xi} =\|J^s_x \varphi \|_{L^2_\lambda} \quad ,
$$
where $ \langle \cdot \rangle = (1+|\cdot|^2)^{1/2} $ and $
\widehat{J^s_x \varphi}(\xi)=\langle \xi \rangle^{s}
\widehat{\varphi}(\xi)
$. \\
For $ s\ge 0 $, the closed subspace of zero mean value functions of $ H^s_\lambda $ will be denoted by $ \dot{H}^s_\lambda $.\\
The Lebesgue spaces $ L^q_\lambda $, $ 1\le q \le \infty $, will be defined as usually by
$$
\|\varphi\|_{L^q_\lambda}=\Bigl( \int_{\R/(2\pi\lambda) \Z} |\varphi(x) |^q \, dx\Bigr)^{1/q}
$$
with the obvious modification for $ q=\infty $. \\
In the same way, for a function $ u(t,x) $ on $ \R\times\R/(2\pi \lambda) \Z $, we define its space-time Fourier transform by
$$
\hat{u}(\tau,\xi)={\cal F}_{t,x}(u)(\tau,\xi)=\int_{\R} \int_{\R/(2\pi\lambda) \Z} e^{-i (\tau
t+ \xi x)} \, u(t,x) \, dx dt , \quad \forall (\tau,\xi) \in
\R\times \lambda^{-1}\Z \quad .
$$
 We define the  Bourgain spaces $ X^{b,s}_\lambda $, $ Z^{b,s}_\lambda $ and $ Y^s_\lambda $ of $ (2\pi\lambda) $-periodic (in $x$) functions respectively
  endowed with the norm
\begin{eqnarray}
\| u \|_{X^{b,s}_{\lambda} }  =
 \| \langle \tau+\xi |\xi|\rangle^{b}  \langle \xi \rangle^s
  \hat{u}\|_{L^2_{\tau,\xi}} =
  \| \langle \tau\rangle^{b}  \langle \xi \rangle^s
  {\cal F}_{t,x}(V(t) u ) \|_{L^2_{\tau,\xi}}\; ,
\end{eqnarray}
\begin{eqnarray}
\| u \|_{Z^{b,s}_{\lambda} } =
 \|  \langle \tau+\xi |\xi|\rangle^b \langle \xi \rangle^s
  \hat{u} \|_{L^2_\xi L^1_\tau}= |  \langle \tau\rangle^b \langle \xi \rangle^s
  {\cal F}_{t,x}(V(t) u ) \|_{L^2_\xi L^1_\tau}\; ,
\end{eqnarray}
and
\begin{equation}
\|u\|_{Y^s_\lambda} = \|u\|_{X^{1/2,s}_\lambda} +
\|u\|_{Z^{0,s}_\lambda} \quad .
\end{equation}
 For $ T > 0$ and a function space $ B_\lambda= X^{b,s}_{\lambda},
  \, Z^{b,s}_{\lambda} $ or $ Y^s_\lambda $,
   we denote by  $ B_{T,\lambda} $ the corresponding restriction in time   space
 endowed with the norm
$$
\| u \|_{B_{T,\lambda}} =\inf_{w\in B_{\lambda}} \{ \|
w\|_{B_{\lambda}} , \, w(\cdot)\equiv u(\cdot) \hbox{ on } [0,T] \, \}\;.
$$
Recall that $ Y^s_{T,\lambda} \hookrightarrow Z^{0,s}_{T,\lambda} \hookrightarrow C([0,T];H^s_\lambda) $. \\
 $ L^p_t L^q_\lambda $ and $ L^p_T L^q_\lambda $ will denote the Lebesgue spaces
$$
\|u\|_{L^p_t L^q_\lambda}=\Bigl( \int_{\R}
\|u(t,\cdot)\|_{L^q_\lambda}^p \, dt \Bigr)^{1/p}\quad \mbox{ and
}\quad \|u\|_{L^p_T L^q_\lambda}=\Bigr( \int_{0}^T
\|u(t,\cdot)\|_{L^q_\lambda}^p \, dt\Bigr)^{1/p}
$$
with the obvious modification for $ p=\infty $. \\
Finally, for all function spaces of $ (2\pi\lambda) $-periodic functions, we will drop the index $ \lambda $ when $ \lambda=1 $. \\
We will denote by $ P_+ $ and $ P_- $ the projection on respectiveley the positive and the negative spatial Fourier modes. Moreover,
for $ a\ge  0 $, we will denote by $ P_a $ and $ P_{>a} $ the projection on respectively the spatial Fourier modes of absolute value equal or less than $ a$ and the  spatial Fourier modes strictly larger than $ a$.
\subsection{Main result}
\begin{theorem}\label{main}
For all $u_0\in H^{s}(\T)$ with $s\ge 1/2 $ and all $ T> 0 $, there
exists a unique solution
$$
u\in  C([0,T]; H^{s}(\T)) \cap X^{1/2,0}_{T,\lambda}
$$
of the Benjamin-Ono equation (BO). \\
Moreover $ u\in C_b(\R,H^{1/2}(\T)) $ and  the map $ u_0\mapsto u
$  is continuous from $H^{s}(\T)$ into $C([0,T],H^{s}(\T))$ and
Lipschitz on every bounded set from  $ \dot{H}^{s}(\T) $ into
$C([0,T],\dot{H}^{s}(\T))$.
\end{theorem}
\begin{rem}
Actually, we  prove that the flow-map is Lipschitz on every
bounded set of  any hyperplan of $ H^{s}(\T) $ with a prescribed
mean value.
\end{rem}
\begin{rem} \label{rem2}
For KdV-like equation,
\begin{equation}
u_t + \partial_x D_x^{2\alpha} u = u u_x , \quad \alpha\ge 0 , \label{R}
\end{equation}
 one can easily prove (see Proposition \ref{KdVg} in the appendix) that
the map $ u_0\mapsto u $ is not uniformly continuous in $ H^s(\T)
$ for $ s>0 $. This has no relation with the order of the
dispersion and is only related to the nonlinear transport
equation $u_t=u u_x $. More precisely, the key point is that if $
u(t,x) $ is a solution of \re{R} then $ u(t,x+\omega t) +\omega
$, where $ \omega $ is any constant,  is also a solution (see
\cite{KT2} or \cite{nico}).
 Therefore, we think that a good notion of stability for this type of equations is the uniform continuity of the flow-map on bounded sets of hyperplans of functions with fixed mean value.
 Indeed, the restriction on such hyperplans prevent to perform  the above transformation. \\
 In this point of view, (BO) shares the same behavior as KdV since by Theorem \ref{main} and \cite{Bo2},  the flow-map for these both equations is Lipschitz on any bounded balls of such hyperplans. This is in sharp contrast with the real line case where the flow-map of KdV is uniformly continuous on every bounded set of $ H^s(\R) $, $s\ge 0 $ (\cite{Bo2}),  whereas the flow-map of (BO) is not uniformly continuous on the unit ball of $ H^s(\R) $ for $ s> 0 $ (\cite{KT2}).
\end{rem}
The main tools to prove Theorem \ref{main} are the gauge
transformation of T. Tao and the Fourier restriction spaces
introduced by Bourgain. Recall that in order to solve (BO), T.
Tao \cite{T} performed a kind of complex Cole-Hopf
transformation\footnote{Note that projecting (BO) on the non
negative frequencies, one gets the following equation : $
\partial_t (P_+ u) -i\partial_x^2 P_+ u =- P_+(u u_x) $} $
W=P_+(e^{-iF/2}) $, where $ F $ is a primitive of $ u $. In the
periodic setting, requiring that $ u $ has  mean value zero, we
can take  $ F=
  \partial_x^{-1} u $ the unique zero mean value primitive of $
  u$. By the mean value theorem, it is then easy to check
  that the above gauge transformation  is Lipschitz from $ L^2_\lambda $ to $
  L^\infty_\lambda $. This property, which is not true on the real line, is crucial to derive the
  regularity of the flow-map. On the other hand, when one expresses  $ u
  $ in terms of $ w=\partial_x W  $ one gets something like
  $$
  u=e^{iF} w +R(u)
  $$
  which is not so good since multiplication by  gauge function as
  $ e^{iF} $ behaves not so well in Bourgain spaces\footnote{Let us note that Bourgain spaces do not
   enjoy an algebra property}. For this reason we will take $ u $ and $ w $ with
    different space
  regularities in the scale of Bourgain spaces. Indeed, whereas $ w $
  will belong to
   $ Y^{1/2}_{T,\lambda}$, $ u $ will belong to $ X^{1/2,0}_{T,\lambda}\cap C([0,T];H^{1/2}(\T))
   $. Note that this bad behavior is compensated by the fact that
   Bourgain spaces permit to gain one derivative when estimating  the
  nonlinear term of  the  equation \re{eqw} satisfied by $ w$. \vspace*{5mm} \\
  {\it This paper is organized as follows:} In the next section
   we recall some linear estimates in Bourgain type spaces.
    In Section 3 we introduce the gauge transform and establish the key nonlinear
     estimates. Section 4 is devoted to the proof of the local existence result
      for small data in $ \dot{H}^{s}(\T) $ and Section 5 to the uniqueness
       and the regularity of the flow-map. Finally, in Section 6 we use dilation
        arguments to extend the result for arbitrary large data and thus prove
         Theorem \ref{main}. The appendix is divided in 3 parts. In the first part
          we present a proof communicated to us by N. Tzvetkov of the crucial linear estimate
           \re{l5} first proved by Bourgain \cite{Bo1}. In the second part we
            present a brief proof of the lack of uniform continuity of the flow-map
             for KdV type equations on the circle (see Remark \ref{rem2}).
              Finally the third part is devoted to the proof of  technical lemmas
               needed in Section 3.
%%%%%%%%%%%%%%%%%%%%%%%%%%%%%%%%%%%%%%%%%%%%%%%%%%%%%%%%%%%%%%%
%%%%%%%%%%%%%%%%%%%%%%%%%%%%%%%%%%%%%%%%%%%%%%%%%%%%%%%%%%%%%%%
\section{Linear Estimates}
One of the main ingredient is the following linear estimate due to
Bourgain \cite{Bo1}. We present in the appendix a shorter proof
of this result.
\begin{equation}
\|v\|_{L^4_{\lambda^2,\lambda}} \lesssim
\|v\|_{X^{3/8,0}_{\lambda^2,\lambda}} \quad . \label{l5}
\end{equation}
This estimate is proved in \cite{Bo1}  for  Bourgain spaces
associated with the Schr\"odinger group and for a period equal to $
1 $. The corresponding estimate  for the Benjamin-Ono group follows by  writting $ v $ as the sum of its  positive and negative frequency parts. Also the estimate for any period $ \lambda\ge 1
$ follows directly from dilation arguments.  The following classical lemma
(cf. \cite{GTV2}) enables to  deduce a localized version in time of
the above estimate with a gain of a small factor of $ T $.
\begin{lem}
For any  $ T> 0 $ and $ 0 \le b \le 1/2 $, it holds
\label{line1}
\begin{equation}
\|v\|_{X^{b,0}_{T,\lambda}} \lesssim \, T^{1/2-b}\,
\|v\|_{X^{1/2,0}_{T,\lambda}} \quad . \label{l6}
\end{equation}
\end{lem}
Combining \re{l5} and \re{l6} we deduce that for $ 3/8\le b \le
1/2 $ and $ 0<T \le \lambda^2$, it holds
\begin{equation}
\|v\|_{L^4_{T,\lambda}} \lesssim \, T^{b-3/8} \, \|v\|_{X_{T,\lambda}^{b,0}}
\quad . \label{l6b}
\end{equation}
 Let us now state some  estimates  for the free group
and the Duhamel operator. Let  $ \psi\in C_0^\infty([-2,2]) $ be
a time function  such that $ 0\le \psi \le 1 $ and $ \psi\equiv 1
$ on $ [-1,1] $. The following linear estimates are well-known (cf. \cite{Bo1}, \cite{G}).
\begin{lem}
\label{line2} For all  $ \varphi \in {H^s_\lambda}$, it holds :
\begin{equation}
\|\psi(t)V(t) \varphi \|_{Y^{s}_{\lambda}} \lesssim \|\varphi
\|_{H^s_\lambda} \quad . \label{l1}
\end{equation}
\end{lem}
\begin{lem}
\label{line3} For all $ G \in X^{-1/2,s}_{\lambda}\cap  Z^{-1,s}_{\lambda} $, it holds
\begin{equation}
\|\psi(t)\int_0^t V(t-t') G(t') \, dt' \|_{Y^{s}_{\lambda}}
\lesssim \| G \|_{X^{-1/2,s}_{\lambda}}+\| G
\|_{Z^{-1,s}_{\lambda}} \quad . \label{l3}
\end{equation}
\end{lem}
Let us recall that  \re{l3} is a  direct consequences of the
following one dimensional (in time) inequalities (cf. \cite{G}):
for any function $ f\in {\cal S}(\R) $, it holds
$$
\|\psi(t) \int_0^t f(t') \, dt' \|_{H^{1/2}_t} \lesssim
\|f\|_{H^{-1/2}_t} + \Bigl\|\frac{{\cal F}_t(f)}{\langle \tau
\rangle}\Bigr\|_{L^1_\tau} $$ and
$$
\Bigl\|{\cal F}_t \Bigl(\psi(t) \int_0^t f(t') \, dt'
\Bigr)\Bigr\|_{L^1_\tau} \lesssim  \Bigl\|\frac{{\cal
F}_t(f)}{\langle \tau \rangle}\Bigr\|_{L^1_\tau} \quad .
$$
\section{Gauge transform and nonlinear estimates}
\subsection{Gauge transform}
 Let $ \lambda\ge 1 $ and
 $ u $ be a smooth $(2\pi \lambda) $-periodic  solution of (BO) with initial data $ u_0 $.
 In the sequel,  we  assume that $u_0 $ has  mean value zero.
  Otherwise we do  the change of unknown :
 \begin{equation}
 v(t,x)=u(t,x-t \into u) -\into u \label{chgtvar}\quad  ,
 \end{equation}
 where $ \into u=P_0(u) =\frac{1}{2\pi \lambda} \int_{\R/(2\pi \lambda)\Z } u $ is the mean value of $ u $.
 Since $ \into u $ is preserved by the flow, it is easy to see that $ v $ satisfies (BO) with $ u_0-\into u_0 $ as initial data. We are thus reduced to the case of zero mean value initial data.
 We  define  $ F=\partial_x^{-1} u $ which is the periodic, zero mean value, primitive of $ u $,
 $$
 \hat {F}(0) =0 \quad \mbox{ and } \widehat{F}(\xi)=\frac{1}{ i\xi} \hat{u}(\xi) , \quad \xi\in \lambda^{-1}\Z^* \quad .
 $$
Following T. Tao \cite{T}, we introduce the gauge transform
\begin{equation}
W=P_+(e^{-iF/2}) \quad . \label{defW}
\end{equation}
Since $ F $ satisfies $$ F_t +{\cal H} F_{xx}=\frac{F_x^2}{2}-\frac{1}{2} \into F_x^2=
  \frac{F_x^2}{2}-\frac{1}{2} P_0(F_x^2) \quad , $$
we  can then check that $ w=W_x=-\frac{i}{2} P_+(e^{-iF/2} F_x) =-\frac{i}{2} P_+(e^{-iF/2}
u)$ satisfies
\begin{eqnarray}
w_t-iw_{xx} & = & -\partial_x P_+\Bigl[  e^{-iF/2}\Bigl(P_-(F_{xx})-\frac{i}{4} P_0(F_x^2)\Bigr)\Bigr] \nonumber \\
 & =  & -\partial_x P_+ \Bigl(W P_-( u_{x})  \Bigr)+ \frac{i}{4} P_0(F_x^2) w \label{eqw} \quad .
\end{eqnarray}
On the other hand, one can write $ u $ as
\begin{equation}
 u =  e^{iF/2} e^{-iF/2} F_x  = 2 i \, e^{iF/2}\partial_x (e^{-iF/2})= 2
i e^{iF/2} w  + 2 i  e^{iF/2} \partial_x P_- (e^{-iF/2} ) \quad , \label{A3}
\end{equation}
thus
\begin{eqnarray}
P_{>1} u & = & 2 i P_{>1} \Bigl( e^{iF/2}w\Bigr) + 2 i P_{>1} \Bigl( e^{iF/2} \partial_x P_-(e^{-iF/2}) \Bigr) \nonumber \\
& =  &  2 i P_{>1} \Bigl( e^{iF/2 }w\Bigr) + 2 i P_{>1} \Bigl(
P_{> 1}(e^{iF/2}) \partial_x  {P_-}(e^{-iF/2} ) \Bigr) \quad . \label{A4}
\end{eqnarray}
The remaining of this section is devoted to the proof of  the following crucial nonlinear estimates on $ u $ and $ w $.
\begin{pro} \label{NON}
Let $  u\in X^{1/2,0}_{T,\lambda}\cap L^\infty_T
\dot{H}^{s}_\lambda $ be a solution of (BO) and $  w\in
Y^{s}_{T,\lambda} $ satisfying \re{eqw}-\re{A3}. Then for $ 0< T
\le 1 $ and $  0\le s\le 1 $, it holds
\begin{equation}
\|w\|_{Y^{s}_{T,\lambda}} \lesssim (1+\|u_0\|_{L^2_\lambda})
\|u_0\|_{H^{s}_\lambda} + T^{1/8}
\|w\|_{X^{1/2,s}_{T,\lambda}}
 \Bigl(\|u\|_{X^{1/2,0}_{T,\lambda}}+  \|u\|_{X^{1/2,0}_{T,\lambda}}^2 \Bigr)\quad . \label{nonlinear1}
 \end{equation}
\begin{eqnarray}
\| u\|_{X^{1/2,s-1/2}_{T,\lambda}} & \lesssim   & T^{1/2}
\|u\|_{L^\infty_T H^{s}_\lambda} +\|u\|_{L^4_{T,\lambda}} \|J_x^s
u\|_{L^4_{T,\lambda}} \quad .
    \label{nonlinear2}
\end{eqnarray}
 Moreover for $ 1/2\le s\le 1 $ and  $ (p,q)= (\infty,2) $ or $ (4,4) $,
\begin{equation}
\|J_x^s u\|_{L^p_T L^q_\lambda} \lesssim \|u_0\|_{L^2_\lambda} +
(1+\|u\|_{L^\infty_T H^{1/2}_\lambda}) \Bigl( \|w\|_{Y^{s}_{T,\lambda}}
+  \|u\|_{L^\infty_T H^{1/2}_\lambda}^2\Bigr)
 \quad .
\label{nonlinear3}
 \end{equation}

\end{pro}
\begin{rem}
It is worth noting  that  \re{nonlinear3}  can also be rewritten
in a convenient
   way for $ s\ge 0 $. We choose the above
 expression involving the $ L^\infty_T H^{1/2}_\lambda $-norm only for simplicity.
The restriction $ s\ge 1/2 $ in our well-posedness result is due to the lost of one half derivative in \re{nonlinear2}
 which can be explained by the bad behavior of Bourgain spaces with respect to multiplication (see \re{A3}).
\end{rem}
\subsection{Proof of Proposition \ref{NON}}
Let us  first prove \re{nonlinear3}. In this purpose, we need the
two following lemmas proven in the appendix (see also \cite{MR2}
and \cite{MR3}). The first one treat the multiplication with the gauge function
 $ e^{-iF/2} $ in Sobolev spaces. The second one shows that, due to the frequency projections, we can share derivatives when taking the $ H^{s} $-norm of the second term of the right-hand side
 to \re{A4}.
\begin{lem}  \label{non1} Let $ 2\le q<\infty $. Let $ f_1 $ and $ f_2 $ be two real-valued functions of $ \L^q_\lambda $
 with mean value zero  and let $g \in L^q_\lambda $ such that $ J_x^\alpha g\in L^q_\lambda  $
 with $ 0\le \alpha\le 1$. Then

\begin{equation}
\Bigl\| J_x^\alpha \Bigl( e^{\mp i \partial_x^{-1} f_1} g  \Bigr)
\Bigr\|_{L^{q}_\lambda} \lesssim \|J_x^\alpha g\|_{L^q_\lambda}
(1+\|f_1\|_{L^q_\lambda}) \quad , \label{non1a}
\end{equation}
and
\begin{eqnarray}
\Bigl\|J_x^\alpha \Bigl( (e^{- i \partial_x^{-1} f_1} - e^{-i
\partial_x^{-1} f_2}) g  \Bigr) \Bigr\|_{L^q_\lambda} &
\lesssim &
\|J_x^\alpha g\|_{L^q_\lambda} \Bigl(\|f_1 -f_2\|_{L^q_\lambda} \nonumber \\
 & & \hspace*{-15mm}+\|e^{-i\partial_x^{-1} f_1}-e^{-i\partial_x^{-1} f_2}\|_{L^\infty_\lambda}
 (1+\|f_1\|_{L^q_\lambda})\Bigr)\; . \label{non1b}
\quad .
\end{eqnarray}
\end{lem}
\begin{lem} \label{non2}
Let $ \alpha\ge 0 $  and $ 1<q<\infty $ then
\begin{equation}
\Bigl\|D^{\alpha}_x P_+\Bigl( f P_- \partial_x  g\Bigr) \Bigr\|_{L^q_\lambda} \lesssim
\| D^{\gamma_1}_x f \|_{L^{q_1}_\lambda } \,
  \| D^{\gamma_2}_x g \|_{L^{q_2}_\lambda} \quad , \label{estlemP+}
\end{equation}
with $ 1<  q_i<\infty $, $  1/q_1+1/q_2=1/q   $ and $
 \left\{ \begin{array}{l} \gamma_1 \ge \alpha , \; \gamma_2\ge 0 \\
 \gamma_1 + \gamma_2=\alpha+1 \end{array}\right. $ .
\end{lem}
%%%%%%
We first  prove \re{nonlinear3}. Since $ u $ is real-valued, it holds
$$
\|J_x^s u\|_{L^p_T L^q_\lambda} \lesssim \|P_1 u\|_{L^p_T L^q_\lambda}+ \|D_x^s P_{>1} u \|_{L^p_T L^q_\lambda} \quad .
$$
From \re{A4},  Lemmas \ref{non1}-\ref{non2}, Sobolev
inequalities and \re{l6b}, we infer that for  $0<T\le 1 $, $1/2\le  s\le 1 $ and $
(p,q) =(\infty,2) $ or $ (4,4)$,
\begin{eqnarray*}
\|D_x^s P_{>1} u\|_{L^p_T L^q_\lambda} &  \lesssim &
(1+\|u\|_{L^\infty_T L^q_\lambda}) \|J_x^s w\|_{L^p_T
L^q_\lambda} \\
 & & +
 \|D_x^{s} P_{>1} e^{iF/2}\|_{L^\infty_T L^{2q}_\lambda}\|
 u\|_{L^p_T L^{2q}_\lambda} \\
&  \lesssim &(1+\|u\|_{L^\infty_T H^{1/2}_\lambda})
\|w\|_{Y^s_{T,\lambda}} \\
 & &  +T^{1/p} \|D_x^{s+1/2-1/(2q)} P_{>1} e^{iF/2}\|_{L^\infty_T
L^2_\lambda}\| u\|_{L^\infty_T H^{1/2}_\lambda} \quad ,
\end{eqnarray*}
 with
 $$
\|D_x^{s+1/2-1/(2q)}P_{>1} e^{iF/2}\|_{L^\infty_T
L^2_\lambda}\lesssim \|u\|_{L^\infty_T H^{1/2}_\lambda}(1+\|u\|_{L^\infty_T L^2_\lambda}) \quad .
 $$
On the other hand, by the Duhamel formulation of the equation,
the unitarity of $ V(t) $ in $ L^2_\lambda $, the continuity of $ \partial_x P_1 $ in
$ L^2_\lambda $ and Sobolev inequalities, we get for $ (p,q)=
(\infty,2)$ or $ (4,4) $,
\begin{equation}
\|P_1 u\|_{L^p_T L^q_\lambda}  \lesssim  \|u_0\|_{L^2_\lambda}+
\|u^2\|_{L^1_T L^2_\lambda}
 \lesssim  \|u_0\|_{L^2_\lambda}+T \|u\|_{L^\infty_T H^{1/2}_\lambda}^2\; .
 \label{estP1u}
\end{equation}
This completes the proof of \re{nonlinear3}. \\
To prove \re{nonlinear2} we start by noticing that
$$
\|u\|_{X^{1,s-1}_{T,\lambda}} \lesssim \|J^{s-1}_x u
\|_{L^2_{T,\lambda}} + \|J_x^{s-1}(u_t+{\cal H} u_{xx})
\|_{L^2_{T,\lambda}}  \quad $$  and thus by the equation and Leibniz
rule for fractional derivatives (cf \cite{KPV4}) we deduce that
\begin{equation}
\|u\|_{X^{1,s-1}_{T,\lambda}} \lesssim T^{1/2} \| u
\|_{L^\infty_T L^2_\lambda} + \|u\|_{L^4_{T,\lambda} } \|J_x^s u
\|_{L^4_{T,\lambda}} \quad . \label{new1}
\end{equation}
Interpolating between \re{new1} and the obvious estimate
 $$
\|u\|_{X^{0,s}_{T,\lambda}}
  \lesssim  T^{1/2} \, \|u\|_{L^\infty_T
 H^s_\lambda} \; ,$$ \re{nonlinear2} follows.

 To prove \re{nonlinear1}, we will first
need to establish two non linear estimates.
These estimates enlight the good behavior in Bourgain spaces of the nonlinear term of \re{eqw}. The frequency projections lead to the smoothing relation \re{C5} which enables somehow to gain one derivative.
In the following lemmas we
will assume that the functions are supported in time in $
[-2T,2T] $.
 Moreover,
 since all the norms appearing in the right-hand side on the inequalities only see the size of the module of the Fourier transform, we can always assume
 that all the functions have non negative Fourier transforms.\\
\begin{lem} \label{nonlin1}
For any $ s\ge 0 $ and $ 0<T\le 1 $,
\begin{equation}
\Bigl\|\partial_x P_+ (W P_-(u_{x}))
\Bigr\|_{X^{-1/2,s}_{\lambda} }\lesssim T^{1/4} \|W_x
\|_{X^{1/2,s}_{\lambda}} \|u \|_{X^{1/2,0}_{\lambda}} \label{C1}\quad .
\end{equation}
\end{lem}
{\it Proof. } As we wrote above, we assume that the functions
have time
support in $ [-2T,2T] $ and  non-negative Fourier transforms. \\
 By duality it thus suffices to show that
\begin{eqnarray}
I & =& \Bigl| \int_{A}  \frac{\langle \xi\rangle^{s} \xi\,  \hat{h}(\tau,\xi) \langle \xi_1\rangle^{-s}
  \xi_1^{-1} \hat{f}(\tau_1,\xi_1) \xi_2 \hat{g}(\tau_2,\xi_2)}{\langle \sigma \rangle^{1/2}\, \langle \sigma_1 \rangle^{1/2}
   \langle \sigma_2 \rangle^{1/2}} \Bigr|\nonumber \\
    & \lesssim & T^{1/4} \, \|h\|_{L^2_{t,\lambda}}\, \|f\|_{L^2_{t,\lambda}}\, \|g\|_{L^2_{t,\lambda}} \label{C2}
\end{eqnarray}
where $ (\tau_2,\xi_2)=(\tau-\tau_1,\xi-\xi_1) $,
$$
\sigma=\sigma(\tau,\xi)=\tau+\xi|\xi|, \quad \sigma_i=\sigma(\tau_i,\xi_i) , \quad i=1,2,
$$
and, due to the frequency projections,  the domain of integration
$ A\subset \R^2\times (\lambda^{-1}\Z)^2 $ is given by
$$
A=\{(\tau,\tau_1, \xi,\xi_1)\in \R^2\times(\lambda^{-1}\Z)^2,
\quad \xi\ge 1/\lambda, \, \xi_1\ge 1/\lambda , \, \xi-\xi_1\le
-1/\lambda \quad \} \quad.
$$
 Note that in the domain of integration above,
\begin{equation}
\xi_1 \ge |\xi_2| \quad \mbox{ and } \quad \xi_1\ge \xi \quad . \label{C3}
\end{equation}
We thus get
\begin{equation}
I\lesssim  \int_{A} \frac{ \xi^{1/2}\,  \hat{h}(\tau,\xi)  \hat{f}(\tau_1,\xi_1) |\xi_2|^{1/2} \hat{g}(\tau_2,\xi_2)}{\langle \sigma \rangle^{1/2}\, \langle \sigma_1 \rangle^{1/2}
   \langle \sigma_2 \rangle^{1/2}}\label{C4}
\end{equation}
Moreover, in $ A $, we have
\begin{equation}
\sigma_1+\sigma_2-\sigma=\xi_1^2-\xi_2^2-\xi^2=-2\xi_2 \, \xi \quad . \label{C5}
\end{equation}
Therefore, if $ \sigma $ is dominant then by Plancherel and \re{l6b},
\begin{eqnarray}
I & \le & \int_{A} \frac{ \hat{h}(\tau,\xi) \,  \hat{f}(\tau_1,\xi_1) \, \hat{g}(\tau_2,\xi_2)}{ \langle \sigma_1 \rangle^{1/2}
   \langle \sigma_2 \rangle^{1/2}} \nonumber \\
    & \lesssim &  \, \|h\|_{L^2_{t,\lambda}}\,  \Bigl\|{\cal F}^{-1}\Bigl( \frac{\hat{f}}{\langle\sigma \rangle^{1/2}}\Bigr)
    \Bigr\|_{L^4_{t,\lambda}}\, \Bigl\|{\cal F}^{-1}\Bigl( \frac{\hat{g}}{\langle\sigma \rangle^{1/2}}\Bigr)\Bigr\|_{L^4_{t,\lambda}} \nonumber \\
     & \lesssim & T^{1/4} \, \|h\|_{L^2_{t,\lambda}}\, \|f\|_{L^2_{t,\lambda}}\, \|g\|_{L^2_{t,\lambda}} \quad .\label{C6}
\end{eqnarray}
Finally, it is clear that the cases $ \sigma_1 $ and $ \sigma_2 $ dominant can be treated in exactly the same way.
\begin{lem} \label{nonlin2}
For any $ s\ge 0 $ and $ 0<T \le 1 $,
\begin{equation}
\Bigl\| \partial_x {P}_+ (W P_-(u_{x}))\Bigr)
\Bigr\|_{Z^{-1,s}_{\lambda}} \lesssim T^{1/8} \, \|W_x
\|_{X^{1/2,s}_{\lambda}} \|u \|_{X^{1/2,0}_{\lambda}} \quad .
\label{C7}
\end{equation}
\end{lem}
{\it Proof.}
 First note that   by Cauchy-Schwarz in $ \tau $,
$$
\Bigl\|\frac{\langle \xi \rangle^{s}}{ \langle \sigma \rangle} {\cal F}
\Bigl( \partial_x {P}_+ (W P_-(u_{x}))\Bigr)  \Bigr\|_{L^2_\xi L^1_\tau} \lesssim
\Bigl\|\partial_x {P}_+ (W P_-(u_{x}))\Bigr)  \Bigr\|_{X^{-1/2+\varepsilon,s}_\lambda}, \; \varepsilon>0 ,
$$
which can be estimated in the regions
$ \{\langle \sigma_1\rangle \ge \langle \sigma \rangle/10\} $ and
$ \{\langle \sigma_2\rangle \ge \langle \sigma \rangle/10\} $  as in the proof of the preceding lemma (with $ T^{1/8} $ instead of $ T^{1/4} $), since we only need the weight $ \langle \sigma \rangle^{-3/8} $
 in those regions. \\
 Moreover, in the region $ \{\xi_1\le 1 \} $, using \re{C3} and then \re{l6b} it is easy to show that
\begin{eqnarray*}
\Bigl\|\partial_x {P}_+ (W P_-(u_{x}))\Bigr)  \Bigr\|_{X^{-1/2+\varepsilon,s}_\lambda} &\lesssim &
 \Bigl\| W_x u \Bigr\|_{L^2_{t,\lambda}} \nonumber \\
   & \lesssim  & \|W_x\|_{L^4_{t,\lambda }}  \|u\|_{L^4_{t,\lambda }}\nonumber \\
     & \lesssim & T^{1/4} \, \|W_x\|_{X^{1/2,0}_\lambda} \, \|u\|_{X^{1/2,0}_\lambda} \quad .
\end{eqnarray*}
%%%%%%%%%
 Now in the region  $ \{ \langle \sigma \rangle \ge  10 \langle \sigma_1 \rangle , \quad
\langle \sigma \rangle \ge 10 \langle \sigma_2 \rangle, \; \xi_1>1 \, \} $, we proceed as in \cite{CKSTT2}.
By \re{C5}, in this region we  have :
\begin{equation}
\langle \sigma \rangle \sim  \langle \xi \, \xi_2 \rangle \label{C9} \quad .
\end{equation}
By symmetry we can moreover assume that $ |\sigma_1| \ge |\sigma_2| $. We note that proving \re{C7} is equivalent to proving
\begin{equation}
I\lesssim \|f\|_{L^2_{\tau,\xi}} \|g\|_{L^2_{\tau,\xi}} \quad, \label{C9bis}
\end{equation}
where
\begin{equation}
I  =
 \Bigl\| \chi_{\{\xi\ge 1/\lambda\} } \int_{B(\tau,\xi)}
\frac{\langle \xi\rangle^{s}  \xi \langle \xi_1\rangle^{-s}
  \xi_1^{-1} \hat{f}(\tau_1,\xi_1) |\xi_2| \hat{g}(\tau_2,\xi_2)}{\langle \sigma\rangle \, \langle \sigma_1 \rangle^{1/2}
   \langle \sigma_2 \rangle^{1/2}}\Bigr\|_{L^2_\xi L^1_\tau} \label{C9ter}
\end{equation}
 with $ B(\tau,\xi)\subset \R\times \lambda^{-1}\Z $ given by
$$
B(\tau,\xi)=\{(\tau_1,\xi_1)\in \R\times \Z/\lambda, \, \xi_1\ge 1, \, \xi-\xi_1\le -1/\lambda ,\,
 \langle \sigma \rangle \ge 10 \langle \sigma_1 \rangle, \, |\sigma_1|\ge |\sigma_2| \} \quad .
$$
Recall that \re{C3} holds in the domain of integration of \re{C9ter}.
We divide this  domain into 2 subregions. \\
$ \bullet $  The subregion $  \max(|\sigma_1|,|\sigma_2|) \ge (\xi |\xi_2|)^{\frac{1}{16}} $.
 We will assume that $ \max(|\sigma_1|,|\sigma_2|)=|\sigma_1| $ since the other case can be treated in exactly the same way.
 Then, by \re{C3} and \re{C9}, we get
\begin{equation}
I \lesssim  \Bigl\|  \int_{B_1(\tau,\xi)}
\frac{ \hat{f}(\tau_1,\xi_1) \hat{g}(\tau_2,\xi_2)}{\langle \sigma \rangle^{1/2+\frac{1}{128}} \langle
 \sigma_1 \rangle^{3/8}
   \langle \sigma_2 \rangle^{1/2}}\Bigr\|_{L^2_\xi L^1_\tau} \label{C10}
\end{equation}
where
$$ B_1(\tau,\xi)=\{(\tau_1,\xi_1)\in B(\tau,\xi), \,  \sigma_1 \ge (\xi |\xi_2|)^{\frac{1}{16}} \} $$
and by applying Cauchy-Schwarz in $ \tau $ we obtain thanks to \re{l6b},
\begin{eqnarray}
I  & \lesssim & \Bigl\| \int_{B_1(\tau,\xi)}
\frac{ \hat{f}(\tau_1,\xi_1) \hat{g}(\tau_2,\xi_2)}{ \langle \sigma_1 \rangle^{3/8}
   \langle \sigma_2 \rangle^{1/2}}\Bigr\|_{L^2_{\xi,\tau}}\nonumber \\
    & \lesssim &\Bigl\|{\cal F}^{-1}\Bigl( \frac{\hat{f}}{\langle\sigma \rangle^{3/8}}\Bigr) \Bigr\|_{L^4_{t,\lambda}}\, \Bigl\|{\cal F}^{-1}\Bigl( \frac{\hat{g}}{\langle\sigma \rangle^{1/2}}\Bigr)\Bigr\|_{L^4_{t,\lambda}} \nonumber \\
     & \lesssim & T^{1/8}  \|f\|_{L^2_{t,\lambda}}\, \|g\|_{L^2_{t,\lambda}} \quad .\label{C11}
\end{eqnarray}
$\bullet $ The subregion $  \max(|\sigma_1|,|\sigma_2|) \le (\xi |\xi_2|)^{\frac{1}{16}} $. Changing the $ \tau, \tau_1 $ integrals in $ \tau_1, \tau_2 $ integrals in \re{C9ter} and using \re{C9}, we infer that
$$
I \lesssim \Bigl\|\int_{C(\xi)} \xi_1^{-1}\int_{\tau_1=-\xi_1^2+O(|\xi \, \xi_2|^{1/16})} \frac{\hat{f}(\tau_1,\xi_1)}{\langle \tau_1+\xi_1^2\rangle^{1/2}}
  \int_{\tau_2=\xi_2^2+O(|\xi \, \xi_2|^{1/16})} \frac{\hat{g}(\tau_2,\xi_2)}{\langle \tau_2-\xi_2^2\rangle^{1/2}}\Bigr\|_{L^2_\xi}
$$
with $ C(\xi)=\{\xi_1\in \lambda^{-1}\Z, \xi_1\ge 1, \,
\xi-\xi_1\le -1/\lambda \, \} $. Applying Cauchy-Schwarz inequality in $
\tau_1 $ and $ \tau_2 $ and recalling that $ \xi_1\ge 1 $ we get
$$
I\lesssim\Bigl\|\chi_{\{\xi\ge 1/\lambda\}} \int_{C(\xi)}  \langle \xi_1\rangle^{-1} (\xi |\xi_2|)^\frac{1}{8} K_1(\xi_1) K_2(\xi_2) \Bigr\|_{L^2_\xi}
$$
where
$$
K_1(\xi)=\Bigl( \int_\tau \frac{\hat{f}(\tau,\xi)^2}{\langle \tau+\xi^2 \rangle} \Bigr)^{1/2}\quad \mbox{ and } \quad
K_2(\xi)=\Bigl( \int_\tau \frac{\hat{g}(\tau,\xi)^2}{\langle \tau-\xi^2 \rangle} \Bigr)^{1/2} \quad .
$$
Therefore, by using \re{C3}, H\"older and then Cauchy-Schwarz inequalities,
\begin{eqnarray}
I & \lesssim & \Bigl\|  \langle\xi\rangle^{-\frac{3}{4}} \int_{\xi_1\in \lambda^{-1}
\Z} K_1(\xi_1) K_2(\xi_2) \Bigr\|_{L^2_\xi} \nonumber \\
 & \lesssim & \Bigl\| \int_{\xi_1\in\lambda^{-1}\Z} K_1(\xi_1) K_2(\xi_2) \Bigr\|_{L^\infty_\xi} \nonumber \\
  & \lesssim & \Bigl(\int_{\xi \in \lambda^{-1}\Z} K_1(\xi)^2 \Bigr)^{1/2}\, \Bigl(\int_{\xi_\in
   \lambda^{-1}\Z} K_2(\xi)^2 \Bigr)^{1/2} \nonumber \\
   &  \lesssim & \|f\|_{X^{-1/2,0}_{\lambda}} \, \|g\|_{X^{-1/2,0}_{\lambda}} \nonumber\\
   & \lesssim & T \, \|f\|_{L^2_{t,\lambda}} \, \|g\|_{L^2_{t,\lambda}} \quad ,\label{C12}
\end{eqnarray}
where we used the dual estimate of \re{l6} in the last step. \vspace{2mm} \\
\subsubsection{ Proof of \re{nonlinear1}.}
To complete the proof of \re{nonlinear1} it remains to treat the second term  of the right-hand side of \re{eqw}.
This term is mainly harmless. Indeed by Cauchy-Schwarz inequality in $ \tau $, Sobolev inequalities in time and Minkowski inequality,
 $$
\| P_0(u^2) w \|_{Z_\lambda^{-1,s}} + \| P_0(u^2) w \|_{X^{-1/2,s}_\lambda}   \lesssim
 \| P_0(u^2) w \|_{X^{-1/2+\varepsilon',s}_\lambda}
 \lesssim \| P_0(u^2) w \|_{L^{1+\varepsilon}_t H^s_\lambda} \quad,
 $$
 for some $ 0<\varepsilon,\varepsilon'<\!\!<1 $.
Assuming that $ u$ and $ w$ are supported in time in $ [-2T,2T] $, by H\"older inequality in time and \re{l6b} we get
 $$
\| P_0(u^2) w \|_{L^{1+\varepsilon}_t H^s_\lambda}  \lesssim  T^{1/8} \|J^s_x w \|_{L^4_{t,\lambda}} \|P_0(u^2) \|_{L^2_t L^4_\lambda}
   \lesssim  T^{1/8} \|w\|_{X^{1/2,s}_\lambda} \|P_0(u^2) \|_{L^2_{t,\lambda}} \quad ,
$$
where we used that  $ \|1\|_{L^4_\lambda} \le  \|1\|_{L^2_\lambda} $ since $ \lambda\ge 1 $.
Hence, the following estimate holds:
\begin{equation}
\| P_0(u^2) w \|_{Z_\lambda^{-1,s}} + \| P_0(u^2) w \|_{X^{-1/2,s}_\lambda}   \lesssim
  T^{1/8} \|w\|_{X^{1/2,s}_\lambda} \|u \|^2_{X^{1/2,0}_\lambda} \quad .\label{termsup}
\end{equation}
Now, by the Duhamel formulation of \re{eqw}, for $ 0<T\le 1 $ and $
-T\le t \le T $, we have
\begin{eqnarray*}
w(t)& = & \psi(t) \Bigl[ V(t)w(0) -\int_0^t V(t-t') \partial_x P_+
\Bigl( P_-(\psi_T \, \tilde{u}_x) \psi_T \, \tilde{W}\Bigr)(t') \,
dt'\Bigr] \\
 & & \hspace*{1cm} +\frac{i}{4} \int_0^t V(t-t')\Bigl( P_0(\psi_T^2 \tilde{u}^2) \psi_T \tilde{W}_x\Bigr)(t')\, dt' \, \Bigr] \quad ,
 \end{eqnarray*}
 where $ \psi_T(\cdot)=\psi(\cdot/T) $, $ {\tilde u} $ is an
 extension of $ u $ satisfying $ \|{\tilde
 u}\|_{X^{1/2,0}_\lambda} \le 2 \|u\|_{X^{1/2,0}_{T,\lambda}} $
 and $ \tilde{W} $ is an extension of $ W $
 satisfying  $ \|{\tilde
 W}_x\|_{X^{1/2,s}_\lambda} \le 2 \|W_x\|_{X^{1/2,s}_{T,\lambda}}
 $. At this stage, it is worth noticing that the multiplication by $ \psi_T $
  is continuous in $ X^{1/2,s}_\lambda$ and  $ Y^s_\lambda $
  with a norm which does not depend on $ T>0 $ or $ s\in \R $.
Therefore,
 combining \ Lemmas \ref{line2}-\ref{line3},
\ref{nonlin1}-\ref{nonlin2} and \re{termsup},
  we infer that for $  s\ge 0 $,
$$
\|w\|_{Y^{s}_{T,\lambda}} \lesssim  \|w(0)\|_{H^{s}} + T^{1/8}
\|w\|_{X^{1/2,s}_{T,\lambda}}
 \Bigl( \|u\|_{X^{1/2,0}_{T,\lambda}} +\|u\|_{X^{1/2,0}_{T,\lambda}}^2 \Bigr) \quad .
$$
This proves \re{nonlinear1} since by Lemma \ref{non1}, for $ 0\le
s \le 1 $,
$$
\|w(0)\|_{H^{s}_\lambda} =\|\partial_x P_+ e^{-i\partial_x^{-1} u_0/2}
\|_{H^{s}_\lambda}\lesssim \|u_0 e^{-i\partial_x^{-1} u_0/2}\|_{H^{s}_
\lambda}\lesssim (1+\|u_0\|_{L^2_\lambda})\|u_0\|_{H^{s}_\lambda}
\quad .
$$
\section{Local existence for small data}
We will now prove the local well-posedness result for small data, the result for arbitrary large data  will
  then follow from scaling arguments.
  More precisely, for some small $ 0<\varepsilon<\!\!<1 $ depending only on  the
implicit constant contained in the above estimates\footnote{In
this stage, it worth recalling that
 these  implicit constants do not depend on the period $ \lambda
 $.},  we
will prove  a  local well-posedness result for initial data
belonging
   to the closed ball
  $ B_{\varepsilon,\lambda} $ of $ \dot{H}^{1/2}_\lambda $ defined by
\begin{eqnarray}
B_{\varepsilon,\lambda} & = & \Bigl\{ \varphi\in
\dot{H}^{1/2}_\lambda, \quad \|\varphi \|_{H^{1/2}_\lambda}
\lesssim \varepsilon^2 \, \Bigr\} \label{o1}\quad ,
\end{eqnarray}
with $ \lambda\ge 1 $.
 \subsection{Uniform estimate}
 Let  $ u_0$ belonging to  $ \dot{H}^\infty_\lambda \cap B_{\varepsilon,\lambda} $.
  We want   first to  show that the emanating solution
   $ u\in C(\R;\dot{H}^\infty_\lambda) $, given by the classical well-posedness results (cf. \cite{ABFS}, \cite{Io}),
    satisfies
\begin{equation}
    \| u\|_{X^{1/2,0}_{1,\lambda}}+\|u\|_{L^\infty_1 H^{1/2}_\lambda} \lesssim \varepsilon^2 \quad \mbox{ and }
    \quad  \| w\|_{Y^{1/2}_{1,\lambda}}\lesssim \varepsilon^2 \quad . \label{o2}
    \end{equation}
     Clearly, since $ u $ satisfies the equation, $ u $ belongs in fact to $ C^\infty(\R;H^\infty_\lambda) $.
     Thus, for any $ 0<T\le 1 $,  $   u $ and $ w  $ belong to $ {Y^{\infty}_{T,\lambda}}$ and from the linear estimates
       we easily deduce that
\begin{eqnarray*}
 \| u\|_{X^{1/2,0}_{T,\lambda}} & \lesssim & \|u_0\|_{L^2_\lambda} + \|\partial_x  (u^2)\|_{L^2_{T,\lambda}} \\
  & \lesssim & \|u_0\|_{H^{1/2}_\lambda} +T^{1/2} \|  u \|_{L^\infty_T H^{1}_\lambda}^2
  \end{eqnarray*}
 Recalling also \re{nonlinear1},  by a continuity argument we can thus assume that
 \begin{equation}
     \| u\|_{X^{1/2,0}_{T,\lambda}} + \|u\|_{L^\infty_T H^{1/2}_\lambda}\lesssim \varepsilon \quad \mbox{ and }  \|  w\|_{Y^{1/2}_{T,\lambda}} \lesssim \varepsilon
      \label{o3}
\end{equation}
for some $ 0<T< 1$.  But \re{nonlinear1} then clearly ensures
that $ \| w\|_{Y^{1/2}_{T,\lambda}} \lesssim \varepsilon^2 $
 and \re{nonlinear3}   ensures that
 $$
  \| J^{1/2}_x  u \|_{L^4_{T,\lambda}}+ \|  u\|_{L^\infty_T H^{1/2}_\lambda} \lesssim \varepsilon^2 \quad .
  $$
 We thus deduce from  \re{nonlinear2} that $ \|u\|_{X^{1/2,0}_{T,\lambda}}
 \lesssim  \varepsilon^2 $ and \re{o2} is proven. It then follows
 from \re{nonlinear1} and \re{nonlinear3} that for $ 1/2\le s\le 1 $,
 \begin{equation}
 \|u\|_{L^\infty_1 H^s_\lambda} +\|w\|_{Y^s_{1,\lambda}}\lesssim
 (1+\|u_0\|_{L^2_\lambda}) \|u_0\|_{H^s_\lambda}\quad .
 \label{ets}
 \end{equation}
 \subsubsection{Local existence}
 Let $ u_0\in  B_{\varepsilon,\lambda}\cap H^s_\lambda $ with $1/2\le s \le 1 $ and
 let  $\{ u_{0}^n \} \subset \dot{H}^{\infty}(\T)\cap B_{\varepsilon,\lambda} $  converging to $ u_0 $ in
 $ H^{s}(\T) $. We denote by $ u_n $ the solution of (BO) emanating from $ u_{0}^n $.
 From standard  existence theorems (see for instance \cite{ABFS}, \cite{Io}),  $ u_n\in C(\R;\dot{H}^\infty_\lambda) $.
  According to
 \re{o2},
   $$ \|u_n \|_{X^{1/2,0}_{1,\lambda}}+ \|u_n\|_{L^\infty_1 H^{1/2}_\lambda}\lesssim
   \varepsilon^{2}$$
   and \re{ets} ensures that $$
  \|u_n\|_{L^\infty_1 H^s_\lambda} \lesssim (1+\|u_0\|_{L^2_\lambda}) \|u_0\|_{H^s_\lambda} $$
   uniformly in $ n $. We can  thus  pass to the limit up to a subsequence. We then obtain the
    existence of a solution $ u\in X^{1/2,0}_{1,\lambda}\cap L^\infty_1 \dot{H}^{s}_\lambda $ to the Benjamin-Ono equation with $ u_0 $ as initial data
(there is no problem to pass to  the limit on the nonlinear term
here).
%%%%%%%%%%%%%%%%%%%%%%%%%%%%%%%%%%%%%%%
\section{ Continuity, uniqueness and  regularity of the flow map for small data solutions }
We   are going to   prove   that the flow-map is Lipschitz from  $ B_{\varepsilon,\lambda}\cap {H}^{s}_\lambda $ to $
X^{1/2,0}_{1,\lambda}\cap L^\infty_1 \dot{H}^{s}_\lambda $. The
continuity of $ t\mapsto u (t) $ in $ H^{s}_\lambda $ will follow
directly . So, let  $ u_1 $ and $ u_2 $ be two solutions of (BO) in $
X^{1/2,0}_{T,\lambda}\cap C([0,T]; \dot{H}^{s}_\lambda) $
associated with initial data $ \varphi_1 $ and $ \varphi_2 $
 in $ B_{\varepsilon,\lambda}\cap {H}^{s}_\lambda $.
 We assume that they satisfy
\begin{equation}
   \|u_i\|_{L^\infty_T H^{1/2}_\lambda}+ \| u_i\|_{X^{1/2,0}_{T,\lambda}} \lesssim \varepsilon^2 , \quad i=1,2\quad . \label{oo2}
    \end{equation}
 for some $ 0<T\le 1 $ and  where $ 0<\varepsilon<\!\! <1$ is taken as
above.\\
We set $ W_i=P_+(e^{-iF_i/2}) $ with $ F_i=\partial_x^{-1} u_i $.
\subsection{Regularity and estimate on $ w_i=\partial_x P_+(e^{-i\partial_x^{-1} u_i/2})$}
\label{subsec} The first step consists in showing that $
w_i=\partial_x P_+(e^{-iF_i/2}) $, $ i=1,2 $, belongs to $
Y^{s}_{T,\lambda} $
 and satisfies \re{nonlinear1} with $ u $ and $ u_0 $  replaced by $ u_i $ and $ \varphi_i $. To simplify the notations, we drop the index $ i $ for a while. Since $ u\in C([0,T];H^{1/2}_\lambda) \cap X^{1/2,0}_{T,\lambda} $ and satisfies (BO), $ u_t\in C([0,T];H^{-3/2}_\lambda) $.
  Therefore $ F\in C([0,T];H^{3/2}_\lambda) \cap C^1([0,T];H^{-1/2}_\lambda) \cap X^{1/2,1}_{T,\lambda}$. The following calculations are thus justified:
  \begin{eqnarray*}
  \partial_t W=\partial_t P_+(e^{-iF/2})& = & -\frac{i}{2}P_+(F_t e^{-iF/2}) \\
  & = & -\frac{i}{2} P_+\Bigl(e^{-iF/2} ( -{\cal H} F_{xx}+F_x^2/2 - P_0(F_x^2)/2) \Bigr)
  \end{eqnarray*}
  and
  $$
\partial_{xx} W=\partial_{xx}P_+(e^{-iF/2})=P_+\Bigl(e^{-iF/2} (-F_x^2/4-iF_{xx}/2 ) \Bigr) \quad .
  $$
It follows that $ W $ satifies at least in a distributional sens,
\begin{equation}
W_t-iW_{xx} =-P_+(e^{-iF/2} (P_- F_{xx}-i P_0(F_x^2)/4 )=-P_+(WP_- F_{xx})+\frac{i}{4} P_0(F_x^2) W  \quad .
\label{eqW}
\end{equation}
Therefore $ w=\partial_x W $ satisfies \re{eqw} and $\tilde{W}=\partial_x^{-1}w=W-\into W $ satisfies :
\begin{equation} \label{eqtildeW}
\left\{
\begin{array}{l}
\tilde{W}_t-i\tilde{W}_{xx} =-P_{>0}(\tilde{W}P_- F_{xx})+\frac{i}{4}P_0(F_x^2) \tilde{W} \\
\tilde{W}(0)=W(0)-\into W(0)
\end{array}
\right.
\end{equation}
Since $ F\in C([0,T];H^{3/2}_\lambda)$, one has $ \tilde{W}\in  C([0,T];H^{3/2}_\lambda)\hookrightarrow X^{0,3/2}_{T,\lambda} $.
 Moreover, using Lemma \ref{non2} one can easily check that the right-hand side member of \re{eqtildeW}$_1 $ belongs to $ C([0,T];L^2_\lambda) $. Therefore, by \re{eqtildeW}, $ \tilde{W}\in X^{1,0}_{T,\lambda} $ and by interpolation we deduce that $ \tilde{W}\in Y^{1/2}_{T,\lambda} $. On the other
  hand, on account of estimate \re{nonlinear1} it is easy to construct by a Picard iterative sheme a zero mean value
  solution of \re{eqtildeW}  (for $ F\in X^{1/2,1}_{T,\lambda} $ given) which belongs to $ Y^{s+1}_{T,\lambda} $ and such that $ w=\tilde{W}_x $ satisfies \re{nonlinear1}. Therefore, if
  we prove the uniqueness of the solution $ \tilde{W} $ of \re{eqtildeW} in $ Y^{1/2}_T \cap C([0,T];\dot{H}^{1/2}_\lambda) $,
  for any fixed $ F\in X^{1/2,1}_{T,\lambda} $, we are done. To prove this uniqueness result, we need the two following lemmas, the proof of which
    are slight modifications of those of Lemmas \ref{nonlin1} and \ref{nonlin2}. As in Section 3, we assume that $ F $ and $ \tilde{W} $ are supported in time in $ [-2T,2T] $.
\begin{lem} \label{nonlin1b}
For any  $ 0<T<1 $ and $ \lambda\ge 1 $,
\begin{equation}
\Bigl\| {P}_{>0} (\tilde{W} P_-(F_{xx}))\Bigr)
\Bigr\|_{X^{-1/2,1/2}_{\lambda} }\lesssim T^{1/4} \lambda^{1/2} \, \|\tilde{W}
\|_{X^{1/2,1/2}_{\lambda}} \|F_x \|_{X^{1/2,0}_{\lambda}}
\label{C1b}
\end{equation}
\end{lem}
{\it Proof. } As  mentioned above, the proof is essentially the same as for Lemma \ref{nonlin1}.
 By duality it thus suffices to show that
\begin{eqnarray}
I & =& \Bigl| \int_{A}  \frac{\langle \xi\rangle^{1/2} \,  \hat{h}(\tau,\xi) \langle \xi_1\rangle^{-1/2}
   \hat{f}(\tau_1,\xi_1) \xi_2 \hat{g}(\tau_2,\xi_2)}{\langle \sigma \rangle^{1/2}\, \langle \sigma_1 \rangle^{1/2}
   \langle \sigma_2 \rangle^{1/2}} \Bigr|\nonumber \\
    & \lesssim & T^{1/4} \lambda^{1/2}\, \|h\|_{L^2_{t,\lambda}}\, \|f\|_{L^2_{t,\lambda}}\, \|g\|_{L^2_{t,\lambda}} \label{C2b}
\end{eqnarray}
with the same set $ A $ as in \re{C2}.
 We divide $ A $ into two regions.\\
 $\bullet $ $\xi\ge |\xi_2| $. Then by \re{C3},
 $$
I\lesssim  \int_{A} \frac{   \hat{h}(\tau,\xi)  \hat{f}(\tau_1,\xi_1) |\xi_2| \hat{g}(\tau_2,\xi_2)}{\langle \sigma \rangle^{1/2}\, \langle \sigma_1 \rangle^{1/2}
   \langle \sigma_2 \rangle^{1/2}}
$$
and the result follows, since by \re{C5}, $ |\sigma_1+\sigma_2-\sigma|\ge 2 |\xi_2|^2 $. \vspace{2mm} \\
$ \bullet $ $ \xi\le |\xi_2| $. Since $ \xi\ge 1/\lambda $ in $ A $, we have
\begin{eqnarray*}
 I & \lesssim & \Bigl| \int_{A}  \frac{\langle \xi\rangle^{1/2} \,  \hat{h}(\tau,\xi)
   \hat{f}(\tau_1,\xi_1) |\xi_2|^{1/2} \hat{g}(\tau_2,\xi_2)}{\langle \sigma \rangle^{1/2}\, \langle \sigma_1 \rangle^{1/2}
   \langle \sigma_2 \rangle^{1/2}} \Bigr| \\
  &  \lesssim & \lambda^{1/2} \Bigl| \int_{A}  \frac{ \xi^{1/2} \,  \hat{h}(\tau,\xi)
   \hat{f}(\tau_1,\xi_1) |\xi_2|^{1/2} \hat{g}(\tau_2,\xi_2)}{\langle \sigma \rangle^{1/2}\, \langle \sigma_1 \rangle^{1/2}
   \langle \sigma_2 \rangle^{1/2}} \Bigr|
\end{eqnarray*}
and the result follows thanks to \re{C5}.
\begin{lem} \label{nonlin2b}
For any  $ 0<T<1 $ and $ \lambda\ge 1 $,
\begin{equation}
\Bigl\|  {P}_{>0}  (\tilde{W} P_-(F_{xx}))
\Bigr\|_{Z^{-1,1/2}_{\lambda}} \lesssim T^{1/8} \lambda^{1/2} \, \|\tilde{W}
\|_{X^{1/2,1/2}_{\lambda}} \|F_x \|_{X^{1/2,0}_{\lambda}} \quad .
\label{C7b}
\end{equation}
\end{lem}
{\it Proof.}
The proof is the same as for Lemma \ref{nonlin2} up to some straightforward modifications similar to the ones  of the preceeding lemma. It will thus be ommitted. \vspace{2mm} \\
On the other hand, proceeding exactly as for the obtention of \re{termsup}, it is easy to see that
\begin{equation}
\| P_0(F_x^2) \tilde{W} \|_{Z_\lambda^{-1,s}} + \| P_0(F_x^2) \tilde{W}\|_{X^{-1/2,s}_\lambda}   \lesssim
  T^{1/8} \|\tilde{W} \|_{X^{1/2,s}_\lambda} \|F_x \|^2_{X^{1/2,0}_\lambda} \quad .\label{termsup2}
\end{equation}
Combining \re{eqtildeW}, \re{termsup2}, Lemmas \ref{nonlin1b}-\ref{nonlin2b} and Lemma \ref{line2}-\ref{line3} and proceeding as in the proof of \re{nonlinear1},
we deduce that for
 $ F\in X^{1/2,1}_{T,\lambda} $ fixed there exists a unique  solution
 $\tilde{W} \in Y^{1/2}_T \cap C([0,T];\dot{H}^{1/2}_\lambda) $    of \re{eqtildeW}  with initial data $ \tilde{W}_0\in \dot{H}^{1/2} $. Therefore, by the arguments given above, we can conclude that
  $ w_i=\partial_x P_+ (e^{-iF_i/2}) $ belongs to $ Y^{s}_{T,\lambda} $ and satisfies \re{nonlinear1}.
   In particular,
   by \re{oo2}, we infer that for $ 0< T \le 1 $ and $ i=1,2 $,
   \begin{equation}
   \|w_i\|_{Y^{s}_{T,\lambda}} \lesssim (1+\|\varphi_i\|_{L^2_\lambda}) \|\varphi_i\|_{H^{s}}
   \label{oo3}
   \end{equation}
   and thus thanks to \re{nonlinear3},
   \begin{equation}\label{oo4}
\|w_i\|_{Y^{1/2}_{T,\lambda}}+\|J_x^{1/2} u \|_{L^4_{T,\lambda}}
\lesssim \varepsilon^2 \quad.
   \end{equation}

\subsection{Lipschitz bound  in $ X^{1/2,0}_{T,\lambda} \cap L^\infty_T H^{s}_\lambda$}
We set
 $$
 z=w_1-w_2=-i P_+(e^{-iF_1/2} u_1)+iP_+ (e^{-iF_2/2}
  u_2)
 $$
 with $ F_i =\partial_x^{-1} u_i $. Obviously, $ z $ satisfies
 \begin{eqnarray}
 z_t-i z_{xx}  & = & -\partial_x P_+ \Bigl[ P_-(\partial_x u_1 -\partial_x u_2) W_1\Bigr]
 -\partial_x P_+ \Bigl[ P_-(\partial_x u_2)( W_1-W_2)\Bigr] \nonumber \\
  & & +\frac{i}{4} \Bigl( P_0(u_1^2) z+ P_0(u_1^2-u_2^2) w_2\Bigr)\; . \label{ol3}
 \end{eqnarray}
 On account of  Lemmas \ref{line2}-\ref{line3}, \ref{nonlin1}-\ref{nonlin2} and \re{termsup}, we  thus infer that, for
 $ 0\le s\le 1 $,
 \begin{eqnarray*}
\|z\|_{Y^{s}_{T,\lambda}} & \lesssim  &
\|z(0)\|_{H^{s}_\lambda}+T^{1/8} \Bigl[\|z
\|_{X^{1/2,s}_{T,\lambda}} \Bigl( \|u_1 \|_{X^{1/2,0}_{T,\lambda}}+\|u_1 \|_{X^{1/2,0}_{T,\lambda}}^2\Bigr) \\
 & & +
 \|u_1-u_2\|_{X^{1/2,0}_{T,\lambda}}\|w_2\|_{X^{1/2,s}_{T,\lambda}}\Bigl(1+\|u_1\|_{X^{1/2,0}_{T,\lambda}}+
\|u_2\|_{X^{1/2,0}_{T,\lambda}}\Bigr)\Bigr] \quad .
 \end{eqnarray*}
 Therefore, thanks to \re{oo2} and \re{oo4} for $ 0<T<1
 $,
 \begin{eqnarray}
 \|z\|_{Y^{s}_{T,\lambda}} & \lesssim  & \Bigl(1+\|\varphi_1\|_{H^s_\lambda}(1+
\lambda^{1/2}) \Bigr)\|\varphi_1-\varphi_2\|_{H^{s}_\lambda}
\nonumber \\
& &  + T^{1/8} \,
 \|w_2\|_{X^{1/2,s}_{T,\lambda}} \|u_1-u_2\|_{X^{1/2,0}_{T,\lambda}} \quad
 ,\label{Lip1}
 \end{eqnarray}
 since, by Lemma \ref{non1}, it can be easily seen that
\begin{eqnarray*}
\|z(0)\|_{H^{s}_\lambda} & \lesssim &
\|\varphi_1-\varphi_2\|_{H^{s}_\lambda}\Bigl(1+\|\varphi_1\|_{H^{s}_\lambda}
 +\|\varphi_2\|_{L^2_\lambda}\Bigr)\nonumber\\
 & & +\|e^{-iF_1(0)}-e^{-iF_2(0)}\|_{L^\infty_\lambda} \|\varphi_1\|_{H^{s}_\lambda}
 (1+\|\varphi_1\|_{L^2_\lambda})
\label{o5}
\end{eqnarray*}
with
\begin{eqnarray*}
\|e^{-iF_1(0)}-e^{-iF_2(0)}\|_{L^\infty_{\lambda}}
   \lesssim  \|\partial_x^{-1} (\varphi_1-\varphi_2)\|_{L^\infty_{\lambda}}
   \lesssim
  \lambda^{1/2} \|\varphi_1-\varphi_2 \|_{L^2_\lambda}\quad . \label{o6}
\end{eqnarray*}
On the other hand, writing the equation satisfied by $ v=u_1-u_2
$, proceeding as in \re{new1} and using \re{oo4} we get, for $
0<T<1 $,
\begin{eqnarray}
\|v\|_{X^{1,-1/2}_{T,\lambda}} &  \lesssim  & T^{1/2}\|v\|_{L^\infty_T H^{-1/2}_\lambda}+
 \|\partial_t v+{\cal H}\partial_x^2 v\|_{L^2_{T} H^{-1/2}_\lambda} \nonumber \\
  & \lesssim & \|v\|_{L^\infty_T L^2_\lambda}+\varepsilon^2 \|J_x^{1/2}v\|_{L^4_{T,\lambda}}
\end{eqnarray}
Interpolating this last inequality with the obvious inequality
$$
\|v\|_{X^{0,1/2}_{T,\lambda}}  \lesssim T^{1/2} \|v\|_{L^\infty_T H^{1/2}_\lambda} \; ,
$$
  it follows that, for $ 0<T<1 $,
\begin{equation}
\|v\|_{X^{1/2,0}_{T,\lambda}}  \lesssim \varepsilon^2 \|J_x^{1/2}v\|_{L^4_{T,\lambda}} + \|v\|_{L^\infty_T H^{1/2}_\lambda} \quad . \label{Lip1bis}
\end{equation}
Now, proceeding as in \re{A3}, we infer that
\arraycolsep2pt
\begin{eqnarray*}
v & =& \partial_x F_1 -\partial_x F_2 \nonumber \\
 & = & 2 i e^{iF_1/2} \Bigl[ z +\partial_x P_- \Bigl( e^{-iF_1/2}  -e^{-iF_2/2}
 \Bigr)\Bigr] +2i ( e^{iF_1/2} -e^{iF_2/2}) \Bigl(w_2
  +  \partial_x P_- (e^{-iF_2/2} \Bigr) \label{o7}
\end{eqnarray*}
and thus
 \begin{eqnarray}
P_{>1} v
 & = & 2i P_{>1}(e^{iF_1/2}z)
 +2i P_{>1}\Bigl[  P_{>1}(e^{iF_1/2})\partial_x P_-\Bigl( e^{-iF_1/2} -e^{-iF_2/2} \Bigr)\Bigr] \nonumber \\
& &  \hspace*{-14mm} +2i  P_{>1}\Bigl[( e^{iF_1/2} -e^{iF_2/2})   w_2\Bigr]
  + 2i P_{>1}\Bigl[ P_{>1}( e^{iF_1/2} -e^{iF_2/2}) \partial_x P_- (e^{-iF_2/2} ) \Bigr] \label{o8} \quad .
\end{eqnarray}
\arraycolsep5pt
Therefore, by Lemmas \ref{non1}-\ref{non2}, \re{oo2} and \re{oo4}, for $ 1/2\le s\le 1 $, $ 0<T<1 $ and $ (p,q)=(\infty,2) $ or $ (4,4) $ we get  as in \re{nonlinear3}
\begin{eqnarray}
\|J^s_xv\|_{L^p_T L^q_\lambda} & \lesssim &
\|P_1\,  v\|_{L^p_T L^q_\lambda}+ (1+\|u_1\|_{L^\infty_T
H^{1/2}_\lambda}) \|z\|_{Y^s_{T,\lambda}} \nonumber \\
 & & +\|u_1\|_{L^\infty_T H^{1/2}_\lambda}\| u_1 e^{-iF_1/2} -u_2
 e^{-iF_2/2}\|_{L^\infty_T H^{1/2}_\lambda}+
 \Bigl\|J^s_x\Bigl((e^{iF_1/2}-e^{iF_2/2}) w_2\Bigr) \Bigr\|_{L^p_T L^q_\lambda} \nonumber \\
 &  &+\| u_1 e^{iF_1/2} -u_2 e^{iF_2/2}\|_{L^\infty_T
 H^{1/2}_\lambda}\|u_2\|_{L^\infty_T H^{1/2}_\lambda} \nonumber \\
 & \lesssim & \|P_1\, v\|_{L^p_T L^q_\lambda}+  \|z\|_{Y^s_{T,\lambda}} \nonumber \\
 & &+ \Bigl(\|v\|_{L^\infty H^{1/2}_\lambda} +\|e^{iF_1/2}-e^{iF_2/2} \|_{L^\infty_{T,\lambda}}\Bigr) \Bigl( \|w_2 \|_{Y^s_{T,\lambda}} +\varepsilon^2\Bigr) \; . \label{Lip2}
\end{eqnarray}
Moreover, since the functions $ u_i $, $ i=1,2$, are
real-valued, by the mean value theorem and Sobolev inequalities,
\begin{eqnarray}
\|e^{iF_1/2}-e^{iF_2/2} \|_{L^\infty_{T,\lambda}}  & \lesssim  &
\|\partial_x^{-1} u_1 - \partial_x^{-1} u_2
\|_{L^\infty_{T,\lambda}} \nonumber \\
 & \lesssim & \|\partial_x^{-1} v
\|_{L^\infty_{T} L^2_\lambda}+\| v \|_{L^\infty_{T} L^2_\lambda}
 \end{eqnarray}
and writing the equation satisfied by $ v$, using the unitarity of
$ V(\cdot) $ in $ L^2_\lambda $ and \re{l6b}, it is easily seen that
\begin{eqnarray}
  \|\partial_x^{-1} v
\|_{L^\infty_{T} L^2_\lambda} & \lesssim &\|\partial_x^{-1}
\varphi_1 - \partial_x^{-1} \varphi_2 \|_{ L^2_\lambda}+T^{1/4}
\|v\|_{X^{1/2,0}_{T,\lambda}}(
\|u_1\|_{X^{1/2,0}_{T,\lambda}}+\|u_2\|_{X^{1/2,0}_{T,\lambda}})\nonumber
\\
 & \lesssim &\lambda  \, \| \varphi_1 -  \varphi_2
\|_{ L^2_\lambda}+T^{1/4}\varepsilon^2
\|v\|_{X^{1/2,0}_{T,\lambda}}\; . \label{Lip3}
\end{eqnarray}
Gathering \re{Lip1}, \re{Lip2}, \re{Lip3} and the obvious
estimate (see \re{estP1u})
\begin{equation}
\|P_1\, v\|_{L^p_T L^q_\lambda} \lesssim
\|\varphi_1-\varphi_2\|_{L^2_\lambda} +T \|v\|_{L^\infty_T
H^{1/2}_\lambda} \Bigl( \|u_1\|_{L^\infty_T
H^{1/2}_\lambda}+\|u_2\|_{L^\infty_T H^{1/2}_\lambda}\Bigr),
\end{equation}
we finally deduce that for $ 0<T<1 $ and $ (p,q)=(\infty,2) $ or
$ (4,4) $,
\begin{eqnarray}
\|J^s_x v\|_{L^p_T L^q_\lambda}  & \lesssim & \Bigl(1+\|\varphi_1\|_{H^s_\lambda} (1+ \lambda^{1/2})\Bigr)\| \varphi_1 -  \varphi_2
\|_{ H^s_\lambda} \nonumber \\
&  & \hspace*{-20mm}+\Bigl( \|v\|_{L^\infty_T H^{1/2}_\lambda}+\|v\|_{X^{1/2,0}_{T,\lambda}}+\lambda  \|\varphi_1-\varphi_2\|_{L^2_\lambda}\Bigr)
  \Bigl( \|w_2 \|_{Y^s_{T,\lambda}} +\varepsilon^2\Bigr) \,. \label{Lip4}
\end{eqnarray}
In particular, taking $ s=1/2 $, we deduce from \re{oo2}, \re{oo4} , \re{Lip1bis}  and \re{Lip4} that
\begin{equation}
\|u_1-u_2\|_{L^\infty_T H^{1/2}_\lambda}+\|J^{1/2}_x (u_1-u_2)\|_{L^4_{T,\lambda}}   \lesssim  (1+\varepsilon^2 \lambda)\| \varphi_1 -  \varphi_2
\|_{ H^{1/2}_\lambda} \; . \label{Lip41}
\end{equation}
It then follows from \re{Lip1bis} that
\begin{equation}
\|u_1-u_2\|_{X^{1/2,0}_{T,\lambda}} \lesssim (1+\varepsilon^2 \lambda)\| \varphi_1 -  \varphi_2
\|_{ H^{1/2}_\lambda} \label{Lip5}
\end{equation}
and from \re{oo3} and \re{Lip4} that
\begin{equation} \label{Lip6}
\|u_1-u_2\|_{L^\infty_T H^s_\lambda}  \lesssim  \Bigl(1+(\|\varphi_1\|_{H^s_\lambda}+\|\varphi_2\|_{H^s_\lambda})(1+\lambda)\Bigr)\| \varphi_1 -  \varphi_2
\|_{ H^s_\lambda} \, ,\quad  1/2\le s\le 1  \quad .
\end{equation}
\subsection{Continuity of the trajectory uniqueness and regularity of the flow-map}
Let $ u_0\in B_{\varepsilon,\lambda} \cap H^s_\lambda $, $
1/2\le s\le 1 $.
 With \re{Lip5}-\re{Lip6} in hand, we observe that the approximative sequence
$ u^n $ constructed for the local existence result is a Cauchy sequence
 in  $ C([0,1]; \dot{H}^{s}_\lambda)\cap X^{1/2,0}_{1,\lambda} $
 since $  \|u_{n}\|_{X^{1/2,0}_{1,\lambda}}+
 \|u_n\|_{L^\infty_1 H^{1/2}}  \lesssim \varepsilon^2 $ and
 $ u_{0,n} $ converges to $ u_0 $ in $ \dot{H}^{s}_\lambda$.
  Hence, $ u $ belongs to
  $ C([0,1]; \dot{H}^{s}_\lambda)\cap X^{1/2,0}_{1,\lambda} $. \\
  Now let $ v $ be another solution emanating from $ u_0 $ belonging to the same class  of regularity as $ u $. By Lebesgue monotone convergence theorem, there exists $ k>0 $ such that $ \|(P_{<-k}+P_{>k}) v\|_{X^{1/2,0}_{T,\lambda}} \lesssim \varepsilon^2 $. On the other hand, using Lemma \ref{line2}-\ref{line3}, it is easy to check that
 $$
\| P_{k} v\|_{X^{1/2,0}_{T,\lambda}}  \lesssim     \|u_0\|_{L^2_\lambda}+ k \|v^2\|_{L^2_{T,\lambda}}
 \lesssim    \|u_0\|_{L^2_\lambda}+ T^{1/4}  k \|v\|_{X^{1/2,0}_{T,\lambda}}^2 \quad .
$$
Therefore, for $ T>0 $ small enough we can require that $ v $ satisfies the smallness condition \re{oo2} and
 thus by  \re{Lip6}, $ u_2\equiv u $ on $ [0,T] $. This proves the uniqueness result for  initial data belonging to
   $ B_{\varepsilon,\lambda}$.
Moreover, \re{Lip6} clearly ensures
 that  the  flow-map is
Lipschitz from  $ B_{\varepsilon,\lambda}\cap H^s_\lambda $ into
  $ C([0,1];\dot{H}^{s}_\lambda) $.
  \section{Proof of Theorem \ref{main}}
We used the dilation symmetry argument to extend the result for arbitrary large data.
First note that if $ u(t,x) $ is a $2\pi $-periodic solution of (BO) on $[0,T] $ with initial data
 $ u_0 $ then $ u_\lambda(t,x)=\lambda^{-1} u(\lambda^{-2}t,\lambda^{-1} x)$ is a $ (2\pi \lambda) $-periodic solution of (BO) on
 $ [0,\lambda^2 T] $
 emanating from $u_{0,\lambda}=\lambda^{-1} u_0(\lambda^{-1} x) $.

Let $ u_0\in \dot{H}^{s}(\T) $ with  $ 1/2\le s\le 1 $.  If
  $ u_0 $ belongs to $ B_{\varepsilon,1} $, we are done.  Otherwise, we set
$$
\lambda=\varepsilon^{-4}  \| u_0\|_{H^{1/2}}^2 \ge 1\quad  \mbox{ so that } \quad \|u_{0,\lambda}\|_{H^{1/2}_\lambda} \le \varepsilon^2 \quad .
$$
Therefore,
 $ u_{0,
\lambda} $ belongs to $ B_{\varepsilon,\lambda}$ and we are
reduce to the case of small initial data. Hence, there exists a
unique local solution $ u_\lambda\in
C([0,1];\dot{H}^{s}_\lambda)\cap X^{1/2,0}_{1,\lambda} $ of
(BO)    emanating from $ u_{0,\lambda} $. This proves the
existence
 and uniqueness in $ C([0,T];\dot{H}^{s}(\T))\cap X^{1/2,0}_T $ of the solution $ u $
  emanating from $ u_0 $ with $ T=T(\|u_0\|_{H^{1/2}})
  $ and $ T\rightarrow +\infty $ as $ \|u_0\|_{H^{1/2}} \rightarrow 0 $.
  The global well-posedness result follows  thus directly since, by combining the conservation laws \re{energy} with Sobolev inequalities, it can be easily checked that $ \|u(t)\|_{H^{1/2}} $ remains uniformly bounded for all $ t >0 $.
   The fact that the flow-map is Lipschitz on every bounded set of $ \dot{H}^{s}(\T) $
    follows as well since $\lambda $
    only depends on $ \|u_0\|_{H^{1/2}}$.\vspace*{3mm}\\

 Note that the change of unknown \re{chgtvar} preserves the continuity of
  the solution and the continuity of the flow-map in $ H^{s}(\T)$. Moreover,
  the Lipschitz property (on bounded sets) of the flow-map is also preserved on
  the hyperplans of $ H^{s}(\T) $ with fixed  mean value.
  \begin{rem}
  We prove Theorem \ref{main}  for $ 1/2\le s\le 1$. Of course, the case
  $ s>1 $ can be treated in the same way up to obvious
  modifications. \\
  On the other hand, we think that by the present approach completed with some new ingredients, one can go down to $ H^s(\T) $, $ s>0 $, or even perhaps to $ L^2(\T) $. It would be then very interesting to know if $ L^2(\T) $ is the limit space for the uniform continuity of the flow-map.
  \end{rem}
\section{Appendix}
\subsection{A shorter proof of Estimate \re{l5}}
We present here a very nice shorter proof communicated to us by  Nikolay Tzvetkov of the crucial inequality \re{l5} (with $\lambda=1 $
 and the Schr\"odinger group $ U(\cdot) $)  proven by Bourgain in \cite{Bo1}.
 The estimate for the Benjamin-Ono group can be deduced simply by writing $ v$ as the sum of
 its positive and negative frequency parts. \\
 Note that if $ v(t,x) $ is $2\pi$-periodic in $ t$ and $ x$ then $ U(t) v $ is also $2\pi $-periodic
 in $ t$ and $ x$.  By density of the $2\pi$-periodic functions in $ L^4([0,2\pi]\times \T) $, it thus suffices to prove the result for such functions. We will thus work in the function space $ X^{b,0}(\T^2) $ endowed with the norm
 $$
 \| v\|_{X^{b,0}}=\Bigr(\sum_{\tau\in \Z} \sum_{n\in \Z} \langle \tau-n^2\rangle^{2b}\,  |\widehat{v}(\tau,n)|^2\Bigl)^{1/2}
 \quad .
$$
Let $ v\in X^{b,0}(\T^2) $, $ b\ge 3/8 $. We introduce  a Littlewood-Paley decomposition of $ v $:
$$
v=\sum_{M-dyadic} v_M
$$
where $ \supp \widehat{v_M} \subset \{ (\tau,n)\in \Z^2, \; \langle \tau-n^2 \rangle \sim M \, \} $. Note that
$$
\|v\|^2_{X^{b,0}} \sim  \sum_{M} M^{2b} \|v_M\|_{L^2_{\tau,n}}^2 \quad .
$$
By the triangle inequality, we have
$$
\|v\|_{L^4_{t,x}}^2 =\|v^2\|_{L^2_{t,x}} = \Bigl\| \sum_{M_1,M_2} v_{M_1} v_{M_2} \Bigr\|_{L^2_{t,x}}
 \lesssim \sum_{M_1,M_2} \Bigl\|  v_{M_1} v_{M_2} \Bigr\|_{L^2_{t,x}}
 \lesssim \sum_{M_1\ge M_2} \Bigl\|  v_{M_1} v_{M_2} \Bigr\|_{L^2_{t,x}} \quad .
 $$
 The proof of \re{l5} is based on the following lemma :
 \begin{lem} \label{lemA}
 \begin{equation}
\Bigl\|  v_{M_1} v_{M_2} \Bigr\|_{L^2_{t,x}} \lesssim \Bigl(M_1 \wedge M_2\Bigr)^{1/2}  \Bigl( M_1\vee M_2\Bigr)^{1/4}
 \|v_{M_1} \|_{L^2_{t,x}}\, \|v_{M_2} \|_{L^2_{t,x}} \quad .
 \end{equation}
 \end{lem}
Indeed, with this lemma in hand, rewritting $ M_1 $ as $ M_1=2^l M_2 $ with $ l\in \N $, we get the following chain of inequalities
\begin{eqnarray*}
\sum_{M_1\ge M_2} \Bigl\|  v_{M_1} v_{M_2} \Bigr\|_{L^2_{t,x}} & \lesssim &
\sum_{l\ge 0} \sum_{M_2} M_2^{1/2} (2^l M_2)^{1/4}
 \|v_{M_2} \|_{L^2_{t,x}} \|v_{2^l  M_2} \|_{L^2_{t,x}} \\
 & \lesssim & \sum_{l\ge 0} \sum_{M_2} M_2^{3/8} \|v_{M_2} \|_{L^2_{t,x}}  (2^l M_2)^{3/8} 2^{-l/8}
  \|v_{2^l  M_2} \|_{L^2_{t,x}} \\
   & \lesssim & \sum_{l\ge 0} 2^{-l/8} \Bigl(\sum_{M_2}  M_2^{3/4} \|v_{M_2} \|_{L^2_{t,x}}^2\Bigr)^{1/2}
\Bigl(\sum_{M_2}(2^l M_2)^{3/4} \|v_{2^l M_2} \|_{L^2_{t,x}}^2\Bigr)^{1/2} \\
 & \lesssim  & \|u\|_{X^{3/8,0}}^2 \quad .
\end{eqnarray*}
It thus remains to prove  Lemma \ref{lemA}. By Cauchy-Schwarz in $(\tau_1,n_1) $ we infer that
\begin{eqnarray*}
\|v_{M_1} v_{M_2} \|_{L^2_{t,x}}^2 & =& \sum_{\tau,n} \Bigl|\sum_{\tau_1,n_1} \widehat{v_{M_1}}(\tau_1,n_1)
\widehat{v_{M_2}}(\tau-\tau_1,n-n_1)\Bigr|^2 \\
& \lesssim  & \sum_{\tau,n} \alpha( \tau,n) \sum_{\tau_1,n_1} \Bigl| \widehat{v_{M_1}}(\tau_1,n_1)
\widehat{v_{M_2}}(\tau-\tau_1,n-n_1)\Bigr|^2 \\
& \lesssim & \sup_{\tau,n} \alpha(\tau,n) \, \|v_{M_1}\|_{L^2_{t,x}}^2 \|v_{M_2}\|_{L^2_{t,x}}^2 \quad ,
\end{eqnarray*}
where
\begin{eqnarray*}
\alpha(\tau,n) & = & \# \{ (\tau_1,n_1) , \; (\tau_1,n_1)\in \supp \widehat{v_{M_1}} \mbox{ and }
(\tau_1,n_1)\in \supp \widehat{v_{M_2}} \, \} \\
 & \lesssim & \# \{ (\tau_1,n_1) , \; \langle  \tau_1-n_1^2\rangle \sim M_1 \mbox{ and }
\langle \tau-\tau_1-(n-n_1)^2\rangle \sim M_2 \, \} \\
 & \lesssim & (M_1\wedge M_2) \, \#\{n_1, \; \langle \tau-n_1^2-(n-n_1)^2\rangle \lesssim M_1+M_2 \, \} \\
  & \lesssim & (M_1\wedge M_2) (M_1+M_2)^{1/2} \\
& \lesssim & (M_1\wedge M_2) (M_1\vee M_2)^{1/2} \quad ,
\end{eqnarray*}
since $ \partial^2_{n_1} (\tau-n_1^2-(n-n_1)^2)=-4 $.
\subsection{On the lack of  uniform continuity for KdV type equations on the circle}
We prove the following proposition :
\begin{pro} \label{KdVg}
Let $ (P_\alpha) $, $ \alpha\ge 0 $, denote the Cauchy problem associated with
\begin{equation}
u_t +D_x^{2\alpha} \partial_x u = u u_x , \quad (t,x)\in \R\times \T \\
\end{equation}
The following assertions hold :
\begin{enumerate}
\item Let  $ \alpha\ge 0 $.  For any    $ s>3/2 $ and any $ t>0 $, the flow-map $ u_0\mapsto u(t) $ is not uniformly
continuous from any ball of $ H^s(\T) $ centered at the origin to $ H^s(\T) $.
\item  If $ \alpha\ge 1/2 $ then for any $ s> 0 $ and any  $ T > 0 $, the map $ u_0\mapsto u $ (if it exists) is not uniformly continuous from any ball of
 $ H^s(\T) $ centered at the origin to $ C([0,T];H^s(\T)) $.
 \item  If $ 0\le \alpha < 1/2 $ then the same result holds for any $ s>1/2 $.
 \end{enumerate}
\end{pro}
\begin{rem}
By compactness methods using energy estimates it is known that
 $ (P_\alpha) $, $ \alpha\ge 0 $, is locally well-posed in $ H^s(\T) $
 for $ s>3/2 $ (cf. \cite{ABFS}). Of course this can be improved as soon as $ \alpha $ is large enough. For instance, for $ \alpha=1/2 $(Benjamin-Ono equation) we proved well-posedness for $ s\ge 1/2 $ whereas for $ \alpha=1 $ (KdV) well-posedness is known for $ s\ge -1/2 $ (\cite{KPV2}, \cite{CKSTT1}).
\end{rem}
{\it Proof. } The proof is  a replay of the proof of
Koch-Tzvetkov \cite{KT2} for the Benjamin-Ono equation
 on the real-line (see also \cite{nico} for other dispersive equations). Actually, it is even much simpler. The key reason is that, in the periodic setting, if $ u(t,\cdot) $ is a solution of \re{KdVg} emanating from $ \varphi $  then $ u(t,\cdot+\omega \, t)+\omega $ is {\it exactly} a solution of \re{KdVg} emanating from $ \varphi(\cdot)+\omega $. \\
 For $ \lambda= 2^n $, we set
 $$
 \varphi_\lambda=\lambda^{-s} \sin (\lambda x) ,
 $$
 so that $\varphi_\lambda $ is a $ 2\pi $-periodic function with norm $ \|\varphi_\lambda \|_{H^s}\sim 1 $.
 We denote by $ U_\alpha(\cdot) $ the free group associated with $ (P_\alpha) $ ,
 $$
 \widehat{U_\alpha(t) \varphi}(n)=e^{-i|n|^{2\alpha} n t }\widehat{\varphi}(n) \; , \quad n\in\Z \quad .
 $$
The proposition will be a direct consequence of the fact that the free evolution of $ \varphi_\lambda $ given by
 $$
 (U_\alpha(t) \varphi_\lambda)(x)=\lambda^{-s} \sin \Bigl( -\lambda^{2\alpha+1}t + \lambda x  \Bigr)
 $$
 is a good first approximation on some time interval of the solution $ u_\lambda $ emanating from $ \varphi_\lambda $. More precisely, we have the following key lemma.
 \begin{lem} \label{approx}
Let $ \alpha\ge 0  $ and $ s>1/2 $, then there exists $ 0<\mu<1 $ such that for $ 0<t\lesssim\lambda^{s-(\frac{3}{2}+0^+)} $ the following equality holds in $ H^s(\T) $ :
\begin{equation}
u_\lambda(t,\cdot)=(U_\alpha(t) \varphi_\lambda)(\cdot) +O(\lambda^{-\mu}) \quad . \label{approx1}
\end{equation}
Moreover, if $ \alpha\ge 1/2 $ and  $ 0<s\le 1/2 $,
then there exist  $ 0<\mu<1 $ such that for $ 0<t\lesssim \lambda^{s-1} $ the following equality holds in $ H^s(\T) $ :
\begin{equation}
u_\lambda(t,\cdot)= (U_\alpha(t) \varphi_\lambda)(\cdot)+O(\lambda^{-\mu}) \quad . \label{approx2}
\end{equation}
\end{lem}
Let us  assume this lemma for a while.  For  $ \lambda\ge 1 $, we choose  a time  $ t_\lambda \in [\lambda^{-1+0^+}, \lambda^{s-(\frac{3}{2}+0^+)}] $ in the case $ s>1/2 $ and a time  $ t_\lambda \in [\lambda^{-1+0^+}, \lambda^{s-1}]  $ in the case
$ \alpha\ge 1/2 $ and $ 0<s<1/2 $.  We then set
$$
\varphi_{\lambda,\omega_i}(x)= \varphi_\lambda(x)+\omega_i \quad \mbox{ and }\quad
u_{\lambda,\omega_i}(t,x)= u_\lambda(t,x+\omega_i t )+\omega_i, \quad i=1,2,
$$
with $ \omega_1=(\lambda \, t_\lambda)^{-1}\pi/2 $ and $ \omega_2=-(\lambda \ t_\lambda) ^{-1}\pi/2 $. \\
Obviously $ \|\varphi_{\lambda,\omega_1}-\varphi_{\lambda,\omega_2}\|_{H^s} =(\lambda\, t_\lambda)^{-1}\pi = O(\lambda^{-0^+})
 $. On the other hand, from the triangular
inequality and Lemma \ref{approx} we deduce that
\begin{eqnarray*}
 \|u_{\lambda,\omega_1}(t_\lambda,\cdot)-u_{\lambda,\omega_2}(t_\lambda,\cdot)\|_{H^s}  & \gtrsim &
\|(U_\alpha(t_\lambda)\varphi_\lambda) (\cdot+\omega_1 t_\lambda)-(U_\alpha(t_\lambda)\varphi_\lambda)
 (\cdot+\omega_2 t_\lambda)\|_{H^s} \\
 & &-
\|u_{\lambda,\omega_1}(t_\lambda,\cdot)-[(U_\alpha(t_\lambda)\varphi_\lambda) (\cdot+\omega_1 t_\lambda )  +\omega_1] \|_{H^s} \\
&  &
-\|u_{\lambda,\omega_2}(t_\lambda,\cdot)-[(U_\alpha(t)\varphi_\lambda) (\cdot+\omega_2 t_\lambda ) +\omega_2 ]\|_{H^s}-|\omega_1-\omega_2| \\
& \gtrsim & \|(U_\alpha(t_\lambda)\varphi_\lambda) (\cdot+\omega_1 t_\lambda)-(U_\alpha(t_\lambda)\varphi_\lambda) (\cdot+\omega_2 t_\lambda)  \|_{H^s}\\
 & & - 2 \|u_{\lambda}(t_\lambda,\cdot)-(U_\alpha(t_\lambda)\varphi_\lambda) (\cdot) \|_{H^s}-|\omega_1-\omega_2| \\
& \gtrsim & 1-O(\lambda^{-0^+}) .
\end{eqnarray*}
Where in the last step we computed
\begin{eqnarray*}
\arraycolsep2pt
(U_\alpha(t_\lambda)\varphi_\lambda) (\cdot+\omega_1 t_\lambda ) & -  & (U_\alpha(t_\lambda)\varphi_\lambda
 (\cdot+\omega_2 t_\lambda) \\
 & =& \lambda^{-s}  \Bigl[\sin\Bigl(-\lambda^{2\alpha+1}t_\lambda+\lambda  x +\pi/2\Bigr)
 -\sin \Bigl(-\lambda^{2\alpha+1}t_\lambda +\lambda x- \pi/2\Bigr)  \Bigr]  \\
 &= & 2\lambda^{-s} \cos\Bigl(-\lambda^{2\alpha+1}t_\lambda + \lambda x\Bigr) \quad .
\end{eqnarray*}
Note that, for $ s>3/ 2$, $ t_\lambda $ can be taken arbitrary large when $ \lambda $ goes to infinity.
This proves the proposition since $ u_{\lambda,\omega_i}  $ is the exact solution of $ (P_\alpha )$ emanating from   $ \varphi_{\lambda,\omega_i} $. \vspace{2mm} \\
{\it Proof of Lemma \ref{approx}}
  Obviously, we have
  \begin{eqnarray*}
\partial_t (U_\alpha(t)\varphi_\lambda) +D_x^{2\alpha} \partial_x (U_\alpha(t)\varphi_\lambda)  -(U_\alpha(t)\varphi_\lambda)  \partial_x (U_\alpha(t)\varphi_\lambda)    =
 -(U_\alpha(t)\varphi_\lambda)  \partial_x (U_\alpha(t)\varphi_\lambda)  \\
   =  -\frac{1}{2} \, \lambda^{1-2s}   sin \Bigl(-2 \lambda^{2\alpha+1}t +2 \lambda x  \Bigr) \quad .
   \end{eqnarray*}
   Calling $ F_\lambda $ the term in right-hand side  of the above equality, we thus get
   $$
    \|F_\lambda(t)\|_{L^2} \lesssim \lambda^{1-2s}
, \quad \forall t\in \R \quad .
$$
Setting now $ v_\lambda =u_\lambda-U_\alpha(t)\varphi_\lambda  $, it is easily checked that $ v_\lambda $ verifies
$$
\partial_t v_\lambda +D_x^{2\alpha} \partial_x v_\lambda= v_\lambda \partial_x v_\lambda
+\partial_x (v_\lambda (U_\alpha(t)\varphi_\lambda) ) +F_\lambda \quad .
$$
 Taking the $ L^2 $-scalar product with $ v_\lambda $, we infer that
 $$
 \frac{1}{2} \frac{d}{dt} \|v_\lambda(t)\|_{L^2} \lesssim
  \|\partial_x (U_\alpha(t)\varphi_\lambda)  \|_{L^\infty} \|v_\lambda (t) \|_{L^2} + \|F_\lambda(t) \|_{L^2} \quad .
  $$
Since $ \|\partial_x (U_\alpha(t)\varphi_\lambda)  \|_{L^\infty} \lesssim \lambda^{1-s} $ we thus deduce from Gronwall lemma that for $ 0<T\lesssim\lambda^{s-1} $,
\begin{equation}
\|v_\lambda\|_{L^\infty_T L^2}\lesssim T \|F_\lambda\|_{L^\infty_T L^2} \lesssim \lambda^{1-2s} T \quad . \label{A1}
\end{equation}
On the other hand, classical energy estimates on solutions to  \re{KdVg}  lead  to
$$
\frac{d}{dt} \|u_\lambda (t)\|_{H^{\frac{3}{2}+}} \lesssim \|u_\lambda (t)\|_{H^{\frac{3}{2}+}}^2 \Longrightarrow
 \|u_\lambda(t)\|_{H^{\frac{3}{2}+}} \le \frac{\|\varphi_\lambda\|_{H^{\frac{3}{2}+}}}{1- c\, t\, \|\varphi_\lambda\|_{H^{\frac{3}{2}+}}  } \quad .
 $$
Thus for $ 0<T\le \|\varphi_\lambda\|_{H^{\frac{3}{2}+}}^{-1}/2 \sim\lambda^{s-(\frac{3}{2}+)} $, one has
$$ \|u_\lambda(t)\|_{H^{\frac{3}{2}+}} \lesssim \|\varphi_\lambda\|_{H^{\frac{3}{2}+}} \lesssim \lambda^{(\frac{3}{2}+)-s}
$$
Combining this estimate with the classical energy estimate
$$
\|u_\lambda(t)\|_{H^r} \lesssim \exp(  C \, t \|u_\lambda\|_{H^{\frac{3}{2}+}}) \, \|\varphi_\lambda \|_{H^r} ,\;
\quad r\ge 0 \,,
$$
we  obtain that for $ 0<t\lesssim \lambda^{s-(\frac{3}{2}+)} $,
$$
\|u_\lambda(t)\|_{H^r}\lesssim  \|\varphi_\lambda \|_{H^r} \lesssim \lambda^{r-s} \quad .
$$
Since it is easy to check that for any $ r\in \R $,  $ \|U_\alpha(t)\varphi_\lambda\|_{H^r}  \lesssim \lambda^{r-s}$, we thus infer that for $ r>3/2 $ and $ 0<T\lesssim  \lambda^{s-(\frac{3}{2}+)} $,
\begin{equation}
\|v_\lambda\|_{L^\infty_T H^r}  \lesssim \lambda^{r-s}  \quad . \label{A2}
\end{equation}
For $ s>1/2 $, interpolating between \re{A1} and \re{A2}, with $ r>s $,  yields \re{approx1}.
For  $ \alpha>1/2 $, we use the conservation of the $ L^2 $-norm and of the energy
$$
E(u)=\frac{1}{2}\int_{\R/2\pi\Z} |D_x^\alpha u|^2 +\frac{1}{6}
\int_{\R/2\pi\Z} u^3
$$
to get after some calculations,
$$
\|u(t)\|_{H^\alpha}\lesssim \|\varphi_\lambda\|_{H^\alpha} \lesssim \lambda^{\alpha-s} , \quad \forall t\in \R\, .
$$
 The result for $ 0<s\le 1/2 $ follows by interpolating this inequality with \re{A1}. Finally, for $ \alpha=1/2 $, i.e. the Benjamin-Ono equation, we interpolate \re{A1} with the estimate given by the next conservation law which controls the $ H^1 $-norm of the solution.

\subsection{Proof of Lemmas \ref{non1} and \ref{non2}}
The proofs of Lemmas \ref{non1} and \ref{non2} can be found in the
appendix of \cite{MR2}  in the context of the $ L^p_x L^q_t $
spaces. We present here short proofs for sake of completeness.
%%%%%%%%%%%%%%%%%%%%%%%%%%%%%%%%%%%%%%%%%%%%%%%%%%%%%
\subsubsection{Proof of Lemma \ref{non1}}
We set  $ F= \partial_x^{-1} f_1 $ which is allowed since $ f_1 $
has zero mean value. We first notice that
\begin{eqnarray}
\Bigl\|D_x \Bigl(  (e^{\mp iF})\,  g\Bigr) \Bigr\|_{L^q_\lambda}
& \lesssim &
 \|D_x e^{\mp i F} \|_{L^q_\lambda} \|g\|_{L^\infty_\lambda} +
  \|e^{\mp i F} \|_{L^\infty_\lambda}  \|D_x g\|_{L^q_\lambda} \nonumber \\
  & \lesssim & \|f_1\|_{L^q_\lambda}\|J_x g\|_{L^q_\lambda}+\|D_x g\|_{L^q_\lambda}\nonumber  \\
   & \lesssim & (1+\|f_1\|_{L^q_\lambda})\|J_x g\|_{L^q_\lambda}  \quad , \label{grf}
\end{eqnarray}
where we used that  $f_1 $ is real-valued.
 Interpolating between \re{grf}
 and the obvious inequality
 $ \|  (e^{\mp iF})\,  g \|_{L^q_\lambda} \lesssim \|g\|_{L^q_\lambda}
 $, \re{non1a} follows.

 Finally, \re{non1b} can be obtained exactly in the same way, using that
$$
\|D_x (e^{-i\partial_x^{-1} f_1} -e^{-i\partial_x^{-1} f_2}
)\|_{L^q_\lambda} \lesssim \|e^{-i\partial_x^{-1} f_1}
-e^{-i\partial_x^{-1} f_2} \|_{L^\infty_\lambda}
 \|f_1\|_{L^q_\lambda}+ \|f_1-f_2\|_{L^q_\lambda}\quad .
 $$
%%%%%%%%%%%%%%%%%%%%%%%%%%%%%%%%%%%%%%%%%%%%%%%%%%%%%%%%ù
%%%%%%%%%%%%%%%%%%%%%%%%%%%%%%%%%%%%%%%%%%%%%%%%%%%%%%%%ù
\subsubsection{Proof of Lemma \ref{non2}}

We
need to introduce a few notations to deal with the
Littlewood-Paley decomposition (cf. \cite{KPV4} or \cite{MR2}).
   We
consider $\Delta_k$ and $S_k$ the two operators respectively
defined for $ (2\pi\lambda)$-periodic functions by
$$
  \Delta_k(f)={\cal{F}}^{-1} \Bigl(\eta(2^{-k}\xi) \hat{f}(\xi)\Bigr)
 \quad \mbox{ and }  \quad
S_k(f)={\cal{F}}^{-1} \Bigl(p(2^{-k}\xi) \hat{f}(\xi)\Bigr)\, ,\,
\xi\in \lambda^{-1}\Z  \;,
$$
where $\eta$ is a smooth non negative function supported in $\{
\xi \; , \; 1/2 \leq |\xi |\leq 2\}$, such that ${\displaystyle
\sum_{-\infty}^{+\infty}\eta(2^{-k}\xi )=1}$ for $\xi\neq 0$, and
where
$${\displaystyle p(\xi )=\sum_{j\le -3} \eta(2^{-j}\xi) =1-\sum_{j\ge -2} \eta(2^{-j} \xi
)}$$
 Note that  $p(0)=1 $ and $ \mbox{Supp } p\subset (-1/4,1/4)
$. We also consider the  operators $\tilde{\Delta}_k$ and $\tilde{S}_k
$
 respectively defined by
$$
\tilde{\Delta}_k(f)={\cal{F}}^{-1} \Bigl(\tilde{\eta}(2^{-k}\xi)
\hat{f}(\xi)\Bigr)
 \quad \mbox{ and }  \quad
\tilde{S}_k(f)={\cal{F}}^{-1} \Bigl(\tilde{p}(2^{-k}\xi)
\hat{f}(\xi)\Bigr), \; \xi\in \lambda^{-1}\Z  \,,
$$
where $\tilde{\eta} $ has the same properties as $\eta $ except
that $\tilde{\eta}$ is supported in $\{ \xi \; , \; 1/8 \leq |\xi
|\leq 8\}$ and that $\tilde{\eta}=1$ on $\{ \xi \; , \; 1/4 \leq
|\xi |\leq 4\}$, and where $ \tilde{p}\in C^\infty_0(\R) $ with $
p(\xi)=1 $ for
 $\xi\in [-100,100] $.
. Clearly this implies that
$$
  \forall j\in \Z \; , \; \tilde{\Delta}_j \circ \Delta_j=\Delta_j \;.
$$
$$
\Delta_k(f) \Delta_{k-j}(g)=\tilde{S}_k \Bigl( \Delta_k(f) \Delta_{k-j}(g)\Bigr), \quad
\forall k\in \Z, \; |j|\le 2 \quad .
$$
 We then define the operators
$$
\begin{array}{rr}
\Delta_k^\lambda(f)={\cal F}^{-1}\Bigl(\eta^\lambda (2^{-k} \xi)
f(\xi) \Bigr),
 &
\tilde{\Delta}_k^\lambda(f)={\cal F}^{-1}\Bigl(\tilde{\eta}^\lambda (2^{-k} \xi) f(\xi) \Bigr) ,\\
S_k^\lambda(f)={\cal F}^{-1}\Bigl(p^\lambda (2^{-k} \xi) f(\xi)
\Bigr), & \tilde{S}_k^\lambda(f)={\cal
F}^{-1}\Bigl(\tilde{p}^\lambda (2^{-k} \xi) f(\xi) \Bigr) \quad ,
\end{array}
$$
where $ \eta^\lambda(\xi)=|\xi|^\lambda \, \eta(\xi) $,
 $ \tilde{\eta}^\lambda(\xi)=|\xi|^\lambda \, \tilde{\eta}(\xi) $,
$ p^\lambda(\xi)=|\xi|^\lambda \, p(\xi) $  and
$\tilde{p}^\lambda(\xi)=|\xi|^\lambda \, \tilde{p}(\xi) $.
 Finally we denote by $ {\bold M} $  the maximal operator. \vspace{3mm}

First, note that the zero-modes of $ f $ and $ g $ are not involved  in the expression $P_+(fP_-\partial_x  g)$.
We are thus allowed to use an homogeneous Littlewood-Paley decomposition of $ f$ and $ g $.
Now, since clearly
 $ P_+ (\Delta_l f \, \Delta_k(P_- \partial_x g))=0 $ as soon as $ l\le k-3 $, one has
$$
P_+ \Bigl( \sum_k  S_k(f) \Delta_k(P_- \partial_x g) \Bigr)=0 \quad .
$$
Therefore, the usual homogeneous Littewood-paley decomposition leads to
\begin{eqnarray*}
D^\alpha_x \Bigl( P_+ f P_- \partial_x  g \Bigr) & = & D^\alpha_x P_+
\Bigl[ \sum_k  \Delta_k(f) S_k(P_- \partial_x g) + \sum_{|j|\le 2}\sum_k
 \Delta_{k-j}(f) \Delta_k(P_- \partial_x g)\Bigr]\\
& = & P_+ \Bigl[ \sum_k D_x^\alpha \tilde{\Delta}_k \Bigl( \Delta_k(f)
S_k(P_-\partial_x g)\Bigr)
 + \sum_{|j|\le 2}\sum_k
 D^\alpha_x \tilde{S}_k \Bigl(\Delta_{k-j}(f) \Delta_k(P_- \partial_x g)\Bigr)
\Bigr] \\
& = & P_+ \Bigl[ \sum_k  \tilde{\Delta}_k^\alpha \Bigl(
\Delta_k^{-\gamma_1} (D^{\gamma_1}_x f)
 S_k^{1-\gamma_2}(P_- D^{\gamma_2}_x g)\Bigr) \Bigr] \\
& & + P_+ \Bigl[ \sum_{|j|\le 2}2^{j\gamma_1} \sum_k
\tilde{S}_k^\alpha \Bigl(\Delta_{k-j}^{-\gamma_1}(D^{\gamma_1}_x f)
\Delta_k^{1-\gamma_2} (P_- D^{\gamma_2}_x g)\Bigr) \Bigr] \\
& =& I+II \quad .
\end{eqnarray*}
By the continuity of $ P_+ $ and $ \tilde{S}_k^\alpha $ in $ L^q_\lambda $, Cauchy-Schwarz inequality and Littlewood-Paley square function theorem, for $ 1<q<\infty $,
$$
\| II \|_{L^q_\lambda} \lesssim \|D_x^{\gamma_1} f\|_{L^{q_1}_\lambda} \|D_x^{\gamma_2} g\|_{L^{q_2}_\lambda} \mbox{ with } 1/q_1+1/q_2=1/q , \quad 1<q_i<\infty \, .
$$
On the other hand, by the Littlewood-Paley square function theorem,
$$
\| I\|_{L^q_\lambda} \lesssim
\Bigl\| {\bold M}(P_- D^{\gamma_2}_x g ) \Bigl( \sum_k
|\Delta_k^{-\gamma_1}
 (D^{\gamma_1}_x f) |^2 \Bigr)^{1/2} \Bigr\|_{L^q_\lambda}
$$
and  H\"older inequality and the continuity of the maximal operator on $ L^p_\lambda $, $1<p\le \infty $, yields the result.
%%%%%%%%%%%%%%%%%
\vskip0.3cm \noindent{\bf Acknowledgments .} The author is very grateful to Nikolay Tzvetkov for fruitful discussions on the problem.

\end{document}